\documentclass[11pt]{amsart}

\usepackage{a4wide}
\usepackage{graphicx,eepic}
\usepackage{amsmath,amsfonts,amssymb}
\usepackage{times}

\usepackage{bm}

\newcommand{\la}{\lambda}
\newcommand{\be}{\beta}

\newcommand{\De}{\Delta}
\newcommand{\al}{\alpha}
\newcommand{\ga}{\gamma}
\newcommand{\e}{\varepsilon}

\newcommand{\om}{\omega}
\newcommand{\si}{\sigma}
\newcommand{\Si}{\Sigma}
\newcommand{\BR}{\mathbb{R}}

\newcommand{\A}{\mathcal{A}}
\newcommand{\U}{\mathcal{U}}

\newcommand{\Hh}{\mathcal{H}}
\renewcommand{\S}{\mathcal S}

\newcommand{\C}{\mathcal C}

\newcommand{\fe}{f_{\bm{\e}}}
\newcommand{\xe}{\mathbf{x}_{\bm\e}}

\newcommand{\p}{\mathbf p}
\renewcommand{\a}{\mathbf a}

\newtheorem{lemma}{Lemma}[section]

\newtheorem{prop}[lemma]{Proposition}
\newtheorem{thm}[lemma]{Theorem}
\newtheorem{cor}[lemma]{Corollary}
\theoremstyle{definition}
\newtheorem{Def}[lemma]{Definition}

\theoremstyle{remark}
\newtheorem{rmk}[lemma]{Remark}

\numberwithin{equation}{section}
\numberwithin{table}{section}

\title[Golden Gaskets] {Golden Gaskets: \\ variations on the
  Sierpi\'nski sieve}
\author{Dave Broomhead}
\author{James Montaldi}
\author{Nikita Sidorov}
\date{\today}
\thanks{NS was supported by the EPSRC grant no GR/R61451/01.}

\begin{document}

\begin{abstract}
  We consider the iterated function systems (IFSs) that consist of
  three general similitudes in the plane with centres at three
  non-collinear points, and with a common contraction factor
  $\la\in(0,1)$.

  As is well known, for $\la=1/2$ the invariant set, $\S_\la$, is a
  fractal called the Sierpi\'nski sieve, and for $\la<1/2$ it is also
  a fractal. Our goal is to study $\S_\la$ for this IFS for
  $1/2<\la<2/3$, i.e., when there are ``overlaps" in $\S_\la$ as well
  as ``holes". In this introductory paper we show that despite the
  overlaps (i.e., the Open Set Condition breaking down completely),
  the attractor can still be a totally self-similar fractal, although
  this happens only for a very special family of algebraic $\la$'s
  (so-called ``multinacci numbers").  We evaluate $\dim_H(\S_\la)$ for
  these special values by showing that $\S_\la$ is essentially the
  attractor for an infinite IFS which does satisfy the Open Set
  Condition.  We also show that the set of points in the attractor
  with a unique ``address'' is self-similar, and compute its
  dimension.

  For ``non-multinacci'' values of $\la$ we show that if $\la$ is
  close to $2/3$, then $\S_\la$ has a nonempty interior and that if
  $\la<1/\sqrt{3}$ then $\S_\la$ has zero Lebesgue measure.  Finally
  we discuss higher-dimensional analogues of the model in question.
\end{abstract}

\maketitle

\section*{Introduction and Summary}
\label{sec:intro}

Iterated function systems are one of the most common tools for
constructing fractals.  Usually, however, a very special class of IFSs
is considered for this purpose, namely, those which satisfy the
\emph{Open Set Condition} (OSC)---see below. We present---apparently
for the first time---a family of simple and natural examples of
fractals that originate from IFSs for which the OSC is violated; that
is, for which substantial overlaps occur.

We consider a family of iterated function systems (IFSs) defined
by taking three planar similitudes $f_i(\mathbf x)=\la \mathbf
x+(1-\la)\p_i$ ($i=0,1,2$), where the scaling factor $\la\in(0,1)$
and the centres $\p_i$ are three non-collinear points in
$\mathbb{R}^2$. Without loss of generality we take the centres to
be at the vertices of an equilateral triangle $\Delta$ (see
Section~\ref{sec:final}). The resulting IFS has a unique compact
invariant set $\S_\la$ (depending on $\la$); by definition
$\S_\la$ satisfies
$$
\S_\la = \bigcup_{i=0}^2 f_i(\S_\la).
$$
More conveniently, $\S_\la$ can be found (or rather approximated)
inductively by iterating the $f_j$.  Let
$$
\Delta_n = \bigcup_{\bm\e\in\Sigma^n}\fe(\Delta),
$$
where $\bm\e=(\e_0,\dots,\e_{n-1})\in\Sigma^n$, and
$\Sigma=\{0,1,2\}$, and
$$
\fe = f_{\e_0}\dots f_{\e_{n-1}}.
$$
Since $f_i(\Delta)\subset\Delta$ it follows that
$\Delta_{n+1}\subset\Delta_n$ and then
$$
\S_\la = \lim_{n\to\infty}\Delta_n = \bigcap_{n=1}^\infty \De_n.
$$
In fact all our figures are produced (using Mathematica) by
drawing $\De_n$ for $n$ suitably large, typically between 7 and
10.

\begin{figure}
\centering \scalebox{0.6} {\includegraphics{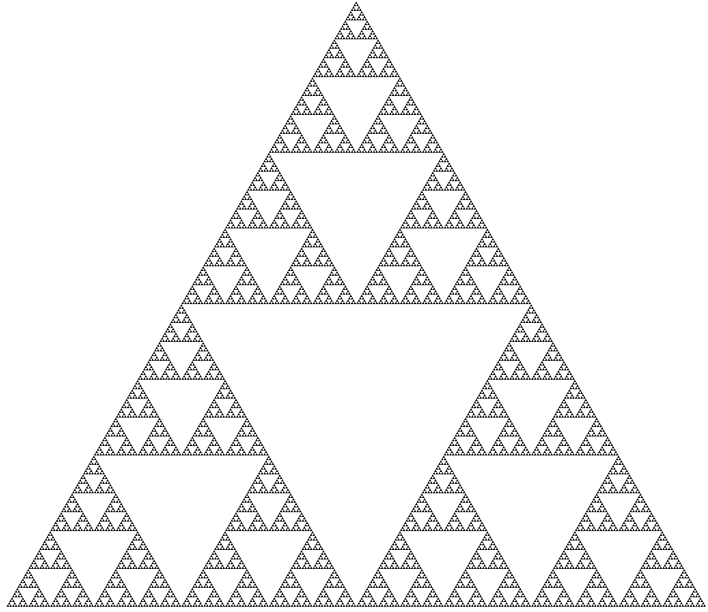}}
\caption{The Sierpi\'nski Sieve.}\label{fig:sg}
\end{figure}

For $\la\le 1/2$ the images of the three similarities are
essentially disjoint (more precisely, the similarities satisfy the
open set condition (OSC)), which makes the invariant set
relatively straightforward to analyse.  For $\la=1/2$ the
invariant set is the famous Sierpi\'nski sieve (or triangle or
gasket)---see Figure~\ref{fig:sg}, and for $\la\le 1/2$ the
invariant set is a self-similar fractal of dimension
$\log3/(-\log\la)$.  On the other hand, if $\la\ge 2/3$ the union
of the three images coincides with the original
triangle\footnote{By ``triangle" we always mean the convex hull of
three points, not just the boundary.} $\De$, so that $\S_\la=\De$.

In this paper we begin a systematic study of the IFS for the
remaining values of $\la$, namely for $\la\in(1/2,\,2/3)$.  In
this region, the three images have significant overlaps, and the
IFS does not satisfy the open set condition, which makes it much
harder to study properties of the invariant set.  For example, it
is not known precisely for which values of $\la$ it has positive
Lebesgue measure.  We do, however,  obtain partial results: for
$\la<1/\sqrt3\approx 0.577$ the invariant set has zero Lebesgue
measure (Proposition~\ref{prop:lebmes}), while for
$\la\ge\la_*\approx 0.648$ it has non-zero Lebesgue measure
(Proposition~\ref{prop:positive}).  We also show
(Proposition~\ref{prop:leb0}) that the Lebesgue measure vanishes
for the specific value $\la=(\sqrt5-1)/2\approx0.618$.

\begin{figure}
\centering \scalebox{0.9}{\includegraphics{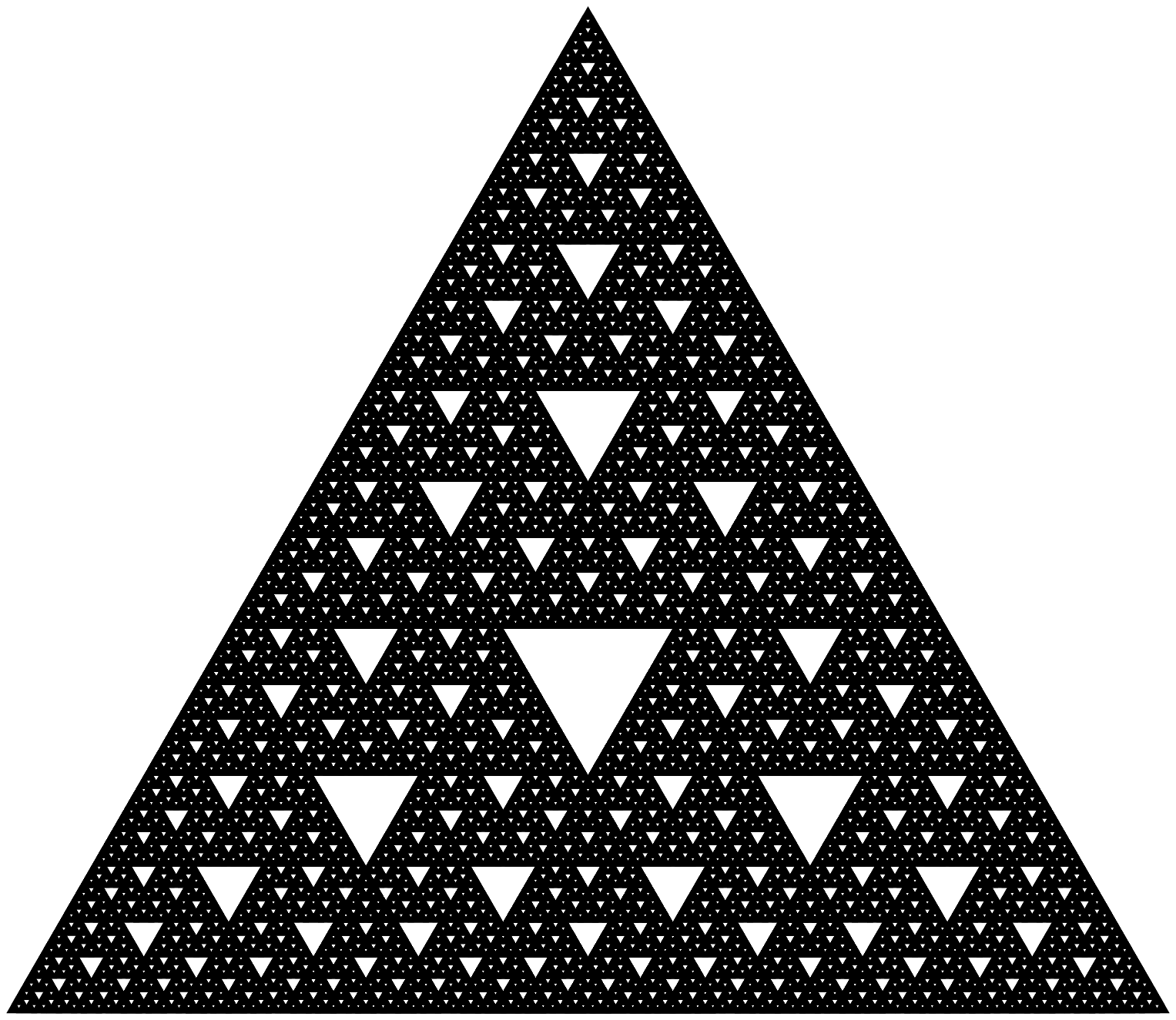}}
\caption{\emph{The} Golden Gasket $\S_{\om_2}$} \label{fig:omega2}
\end{figure}

\begin{figure}
\centering \scalebox{1} {\includegraphics{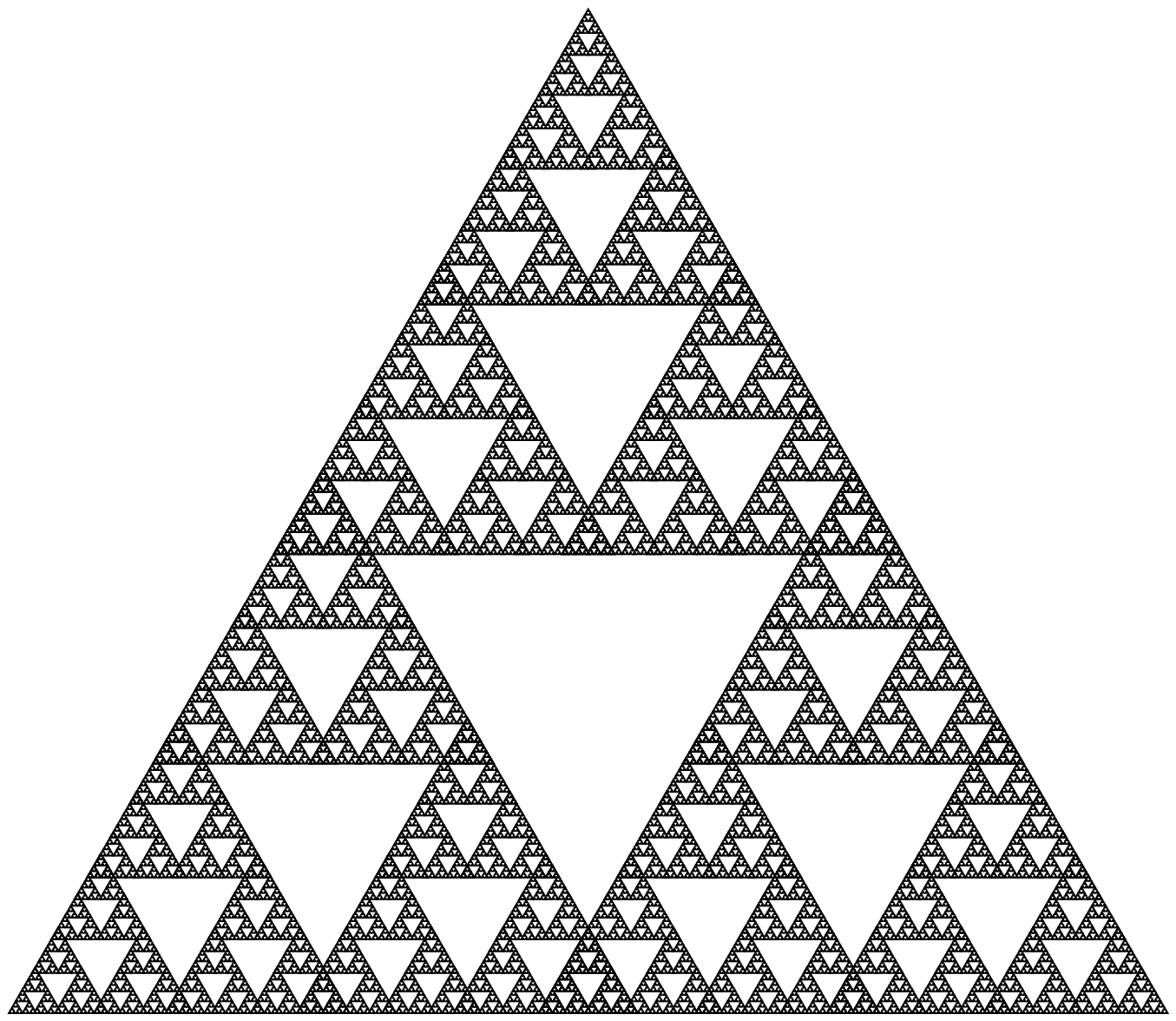}}
\caption{The invariant set $\S_{\om_3}$} \label{fig:omega3}
\end{figure}

\begin{figure}
\centering \scalebox{0.8}{\includegraphics{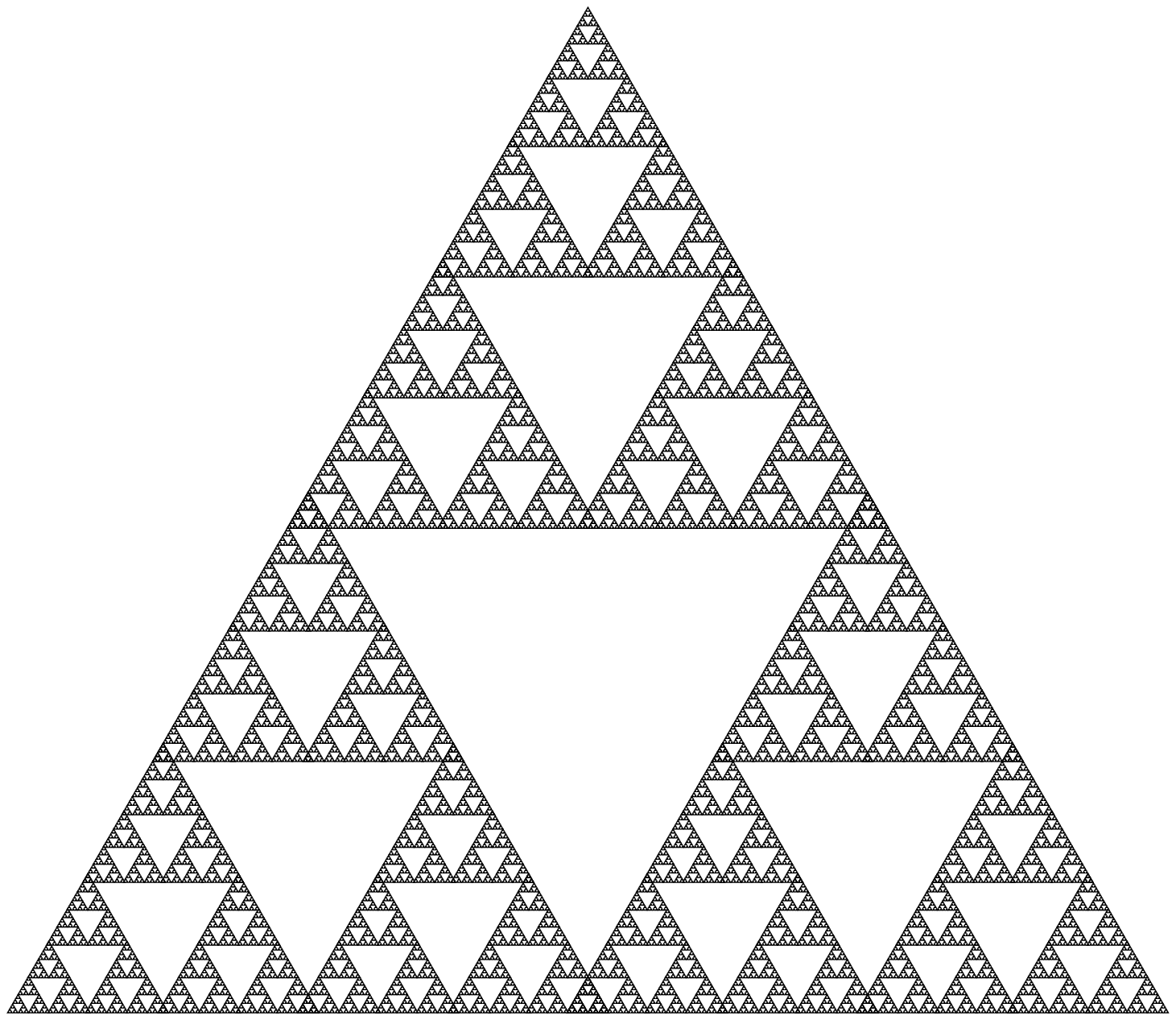}}
\caption{The invariant set $\S_{\om_4}$. Notice the close
resemblance to the Sierpi\'nski sieve in Figure~\ref{fig:sg}.}
\label{fig:omega4}
\end{figure}

The main result of this paper is that there is a countable family
of values of $\la$ in the interval $(1/2,2/3)$---the so-called
\emph{multinacci numbers} $\om_m$---for which the invariant set
$\S_\la$ is totally self-similar (Definition \ref{def:tss}).  We
call the resulting invariant sets \emph{Golden Gaskets}, and the
first three golden gaskets are shown in Figures~\ref{fig:omega2},
\ref{fig:omega3} and \ref{fig:omega4} respectively.  For these
values of $\la$ we are able to compute the Hausdorff dimension of
$\S_\la$ (Theorem~\ref{thm-hd}).  The multinacci numbers are
defined as follows. For each $m\ge2$ the multinacci number $\om_m$
is defined to be the positive solution of the equation
$$
x^m+x^{m-1}+\cdots+x = 1.
$$
The first multinacci number is the golden ratio $\om_2 =
(\sqrt5-1)/2\approx0.618$, and the second is $\om_3\approx0.544$. It
is easy to see that as $m$ increases, so $\om_m$ decreases
monotonically, converging to 1/2.

The key property responsible for the invariant set being totally
self-similar for the multinacci numbers is that for these values
of $\la$ the overlap $f_i(\De)\cap f_j(\De)$ is an image of $\De$,
namely it coincides with $f_if_j^m(\De)$.  On the other hand, we
also show that if $\la\in(1/2,2/3)$ is not a multinacci number
then the invariant set is not totally self-similar
(Theorem~\ref{thm:converse}).

The paper is organized as follows.  In Section~\ref{sec:IFS} we
define the IFS, and introduce the barycentric coordinates we use
for all calculations.  In Section~\ref{sec:holes} we describe the
distribution of holes in the invariant set, and deduce that for
$\la\ge\la_*\approx 0.6478$ the invariant set has nonempty
interior (Proposition~\ref{prop:positive}).

In Section~\ref{sec:gg} we describe explicitly the new family of
\emph{golden gaskets}.  The main result is that for these values
of $\la$ the invariant set is totally self-similar
(Theorem~\ref{thm1}). In Section~\ref{sec:dims} we give several
results on the Lebesgue measure and the Hausdorff dimension of the
invariant set, as described above.  The main result of
Section~\ref{sec:reverse} is that if $\la$ is not multinacci, then
the invariant set is not totally self similar. In
Corollary~\ref{cor:estim} we show how this theorem can be used to
prove a result in number theory---an upper bound for the
``separation constant" that is already known but our proof is very
different, and simpler.

There are two ways to generalize this model: one is to introduce
more similitudes in the plane, and the other is to pass into
higher dimensions, but remaining with simplices (generalizing the
equilateral triangle to higher dimensions).  The first is very
much harder than the second, and in Section~\ref{sec:higher} we
consider the second by way of a very brief discussion of the
``golden sponges'' and a list of a few results that can be
obtained by the same arguments as for the planar case. Finally, in
Section~\ref{sec:final} we end with a few remarks and open
questions.

The appendix contains a detailed proof of the dimensions formula
for the golden gaskets. In doing this, we need to consider points
of the invariant set as determined by a symbolic address: to
$\bm\e\in\Si^\infty$ one associates $\xe\in \S_\la$ by
$\xe=\lim_{n}f_{\e_0}\dots f_{\e_n}(x_0)$ (independently of
$x_0$). The \emph{set of uniqueness} $\U_\la$ consists of those
points in the invariant set that have only one symbolic address in
$\Si^\infty$. For $\la=\om_m$ a multinacci number, we show
$\U_{\om_m}$ to be a self-similar set, and compute its Hausdorff
dimension in Theorem~\ref{thm:dimU}. We also show that ``almost
every" point of $\S_{\om_m}$ (in the sense of prevailing
dimension) has a continuum of different ``addresses"
(Proposition~\ref{prop:contin}).

This is apparently the first paper, where a family of IFSs in
$\BR^2$ with both holes and overlaps is considered in detail. In
$\BR$, however, there has been an attempt to do this, namely, the
famous ``0,1,3"-problem. More precisely, the maps for that model
are as follows: $g_j(x)=\la x+(1-\la)j$, where $x\in\BR$ and
$j\in\{0,1,3\}$. Unfortunately, the problem of describing the
invariant set for this IFS with $\la\in(1/3,2/5)$ (which is
exactly the ``interesting" region) has proved to be very
complicated, and only partial results have been obtained so
far---see \cite{KSS, PS} for more detail.

\section{The iterated function system}
\label{sec:IFS}

Our set-up is as follows. Let $\p_0,\p_1$ and $\p_2$ be the
vertices of the equilateral triangle $\De$:
\[
\p_k=\frac23(\cos(2\pi k/3),\,\sin(2\pi k/3)),\,\, k=0,1,2
\]
(this choice of the scaling will become clear later). Let
$f_0,f_1,f_2$ be three contractions defined as
\begin{equation}\label{eq:ifs0}
f_i(\mathbf x)=\la\mathbf x+(1-\la)\p_i,\quad i=0,1,2.
\end{equation}
Under composition, these functions generate an \emph{iterated function
  system} (IFS)\footnote{Often, in the literature, the term ``IFS"
  means a random functions system endowed with probabilities. Our
  model however will be purely topological.}.

The invariant set (or attractor) of this IFS is defined to be the
unique non-empty compact set $\S_\la$ satisfying
\[
\S_\la=\bigcup_{i=0}^2 f_i(\S_\la).
\]
An iterative procedure exists as follows (see, e.g., \cite{Falc}):
let $\De_0:=\De$ and
\begin{equation}
\De_n:= \bigcup_{i=0}^2 f_i(\De_{n-1}), \quad
n\ge1.\label{eq:deltan}
\end{equation}
The invariant set is then:
\[
\S_\la=\bigcap_{n=0}^\infty\De_n=\lim_{n\to+\infty}\De_n,
\]
where the limit is taken in the Hausdorff metric.

From here on $\Si:=\{0,1,2\}$, $\bm\e$ denotes $(\e_0\dots\e_{n-1})$
(for some $n$) and
\[
\fe := f_{\e_0}\dots f_{\e_{n-1}}.
\]
As is easy to see by induction,
\[
\De_n=\bigcup_{\bm\e\in\Si^n}\fe(\De),
\]
whence $\De_n\subset\De_{n-1}$.


A well studied case is $\la=\frac12$, which leads to the
\emph{Sierpi\'nski sieve} (or \emph{Sierpi\'nski gasket} or
\emph{triangle}) $\S:=\S_{1/2}$---see Figure~\ref{fig:sg}.  Figures
\ref{fig:omega2}, \ref{fig:omega3} and \ref{fig:omega4} show the first
three of the new sequence of fractals, for $\la=\om_2,\om_3$ and
$\om_4$ respectively (the first three multinacci numbers).

\begin{Def}
Recall that the Open Set Condition (OSC) is defined as follows:
let $O$ be the interior of $\De$; then $\bigcup_i f_i(O)\subset
O$, the union being disjoint.
\end{Def}

Note that for $\la=1/2$ the intersections $f_i(\De)\cap f_j(\De)$
($i,j=0,1,2,\ i\neq j$) are two-point sets, i.e., by definition,
this IFS satisfies the Open Set Condition. We would like to
emphasize one more important property of the Sierpi\'nski sieve.
Looking at Figure~\ref{fig:sg}, one immediately sees that each
smaller triangle has the same structure of holes as the big one.
In other words,
\begin{equation}\label{eq:tss}
\fe(\S)=\fe(\De) \cap \S\quad \textrm{for any $\bm\e\in\Si^n$ and
any $n$} .
\end{equation}
\begin{Def}\label{def:tss} We call any set $\S$ that satisfies
  (\ref{eq:tss}), \emph{totally self-similar}.
\end{Def}

Total self-similarity in the case of the Sierpi\'nski sieve
implies, in particular, its holes being well structured: the
$n^{\textrm{th}}$ ``layer" of holes---i.e.,
$\De_n\setminus\De_{n+1}$---contains $3^n$ holes (the central hole
being layer zero), and each of these is surrounded (at a distance
depending on $n$ only) by exactly three holes of the
$(n+1)^{\textrm{th}}$ layer, each smaller in size by the factor
$\la$ ($=1/2$ in this case). Later we will see that only very
special values of $\la$ yield this property of $\S_\la$.

If $\la<1/2$, we have the OSC as well (the intersections
$f_i(\De)\cap f_j(\De),\ i\neq j$ are clearly empty). However, if
$\la\in(1/2,1)$, then $f_i(\De)\cap f_j(\De)$ is always a
triangle, which means that the OSC is not satisfied. This changes
the invariant set dramatically. Our goal is to show that there
exists a countable family of parameters between $1/2$ and 1 which,
despite the lack of the OSC, provide total self-similarity of
$\S_\la$ and, conversely, that for all other $\la$'s there cannot
be total self-similarity.

For technical purposes we introduce a system of coordinates in $\De$
that is more convenient than the usual Cartesian coordinates. Namely,
we identify each point $\mathbf x\in\De$ with a triple $(x,y,z)$,
where
\[
x=\mbox{dist}\,(\mathbf x,[\p_1,\p_2]),\ y=\mbox{dist}\,(\mathbf
x,[\p_0,\p_2]),\ z=\mbox{dist}\,(\mathbf x,[\p_0,\p_1]),
\]
where $[\p_i,\p_j]$ is the edge containing $\p_i$ and $\p_j$.  As
used to be well known from high-school geometry, $x+y+z$ equals
the tripled radius of the inscribed circle, i.e., in our case, 1
(this is why we have chosen the radius of the circumcircle for our
triangle to be equal to $2/3$). These are usually called
\emph{barycentric coordinates}.

\begin{lemma}\label{lem:matr}
  In barycentric coordinates $f_0,f_1,f_2$ act as linear maps.  More
  precisely,
\[
f_0=
\begin{pmatrix}
1 & 1-\la & 1-\la \\ 0 & \la & 0 \\ 0 & 0 & \la
\end{pmatrix},\
f_1=
\begin{pmatrix}
\la & 0 & 0 \\ 1-\la & 1 & 1-\la \\ 0 & 0 & \la
\end{pmatrix},\
f_2=
\begin{pmatrix}
\la & 0 & 0 \\ 0& \la & 0 \\ 1-\la & 1-\la & 1
\end{pmatrix}.
\]
\end{lemma}

\begin{proof}The fact that the $f_i$ are linear in  barycentric
coordinates is a trivial consequence of the $f_i$ being affine.
Let us show that the matrix for $f_0$, say, has the given form
(the proof for $f_1,f_2$ is exactly the same).

Note first that since any vector $(x,y,z)$ is stochastic, so must
be $f_0=(a_{ij})_{i,j=1}^3$, i.e., $a_{ij}\ge0,\ \sum_j a_{ij}=1$
for any $i=1,2,3$. Now, $f_0$ as a map acts as a contraction in
the direction of $\p_0$, which implies that on $y$ or $z$ it acts
simply by multiplying it by $\la$. Therefore,
$a_{21}=0,a_{22}=\la,a_{23}=0,a_{31}=a_{32}=0,a_{33}=\la$, and by
the stochasticity of $f_0$, $a_{11}=1,a_{12}=a_{13}=1-\la$.
\end{proof}

From here on by a \emph{hole} we mean a connected component in
$\De\setminus\S_\la$. First of all, we show that if $\la\ge2/3$, then
there are no holes at all:

\begin{lemma} If $\la\in[2/3,1)$, then $\S_\la=\De$.
\end{lemma}

\begin{proof} It suffices to show that $\bigcup_i f_i(\De)=\De$.
  In barycentric coordinates, $f_0(\De)=\{x\ge1-\la\},
  f_1(\De)=\{y\ge1-\la\}, f_2(\De)=\{z\ge1-\la$\}. For $(x,y,z)$
  to lie in the hole, therefore, the conditions $x<1-\la$,
  $y<1-\la$ and $z<1-\la$ must be satisfied simultaneously. Since
  $\la\ge2/3$ and $x+y+z=1$, this is impossible.
\end{proof}

\section{Structure of the holes}
\label{sec:holes}

Thus, the ``interesting" region is $\la\in(1/2,2/3)$. Let $H_0$
denote the \emph{central hole}, i.e., $H_0=\De\setminus\De_1$; it
is an ``inverted'' equilateral triangle.

\begin{lemma}\label{lem:aux1}
Each hole is a subset of $\bigcup\limits_{\bm\e\in\Si^n}\fe(H_0)$
for some $n\ge1$.
\end{lemma}

\begin{proof} If $\mathbf x$ is in a hole, then there exists $n\ge1$
such that $\mathbf x\in\De_n\setminus\De_{n+1}$. Now our claim
follows from
\[
\De_n\setminus\De_{n+1} =\bigcup_{\bm\e}\fe(\De)\setminus
\bigcup_{\bm\e} \fe(\De_1)
=\bigcup_{\bm\e}\fe(\De\setminus\De_1)=\bigcup_{\bm\e}\fe(H_0).
\]
\end{proof}

\begin{rmk}
  Note that although the $\fe(H_0)$ may not be disjoint, any hole is
  in fact an inverted triangle and a subset of just one of the
  $\fe(H_0)$. We leave this claim without proof, as it is not needed.
\end{rmk}

Let us now derive the formula for any finite combination of $f_i$.
Put
\[
a_k=\begin{cases}1,& \e_k=0\\ 0,& \mathrm{otherwise}\end{cases},\
b_k=\begin{cases}1,& \e_k=1\\ 0,& \mathrm{otherwise}\end{cases},\
c_k=\begin{cases}1,& \e_k=2\\ 0,& \mathrm{otherwise}\end{cases}.
\]
Thus, $a_k,b_k,c_k$ are 0's and 1's and $a_k+b_k+c_k=1$.

\begin{lemma}\label{lem:image1}
Let $\e_k\in\Si$ for $k=0,1,\dots,n$. Then
\[
\fe= \left(
\begin{array}{lll}
(1-\la)\sum\limits_0^{n-1}a_k\la^k + \la^n&
(1-\la)\sum\limits_0^{n-1}a_k\la^k &
(1-\la)\sum\limits_0^{n-1}a_k\la^k \\
(1-\la)\sum\limits_0^{n-1}b_k\la^k &
(1-\la)\sum\limits_0^{n-1}b_k\la^k+\la^n&
(1-\la)\sum\limits_0^{n-1}b_k\la^k \\
(1-\la)\sum\limits_0^{n-1}c_k\la^k &
(1-\la)\sum\limits_0^{n-1}c_k\la^k &
(1-\la)\sum\limits_0^{n-1}c_k\la^k+\la^n
\end{array}\right).
\]
\end{lemma}
\begin{proof} Induction: for $n=1$ this is obviously true; assume
  that the formula is valid for some $n$ and verify its validity for
  $n+1$. Within this proof, we write
  $$
  p_n=(1-\la)\sum_0^{n-1}a_k\la^k,\quad
  q_n=(1-\la)\sum_0^{n-1}b_k\la^k, r_n=(1-\la)\sum_0^{n-1}c_k\la^k.
  $$
  Then by our assumption,
  \begin{align*}
    \fe f_0 &= \left(
      \begin{array}{lll}p_n+\la^n & p_n & p_n\\ q_n & q_n+\la^n & q_n \\
        r_n & r_n & r_n+\la^n
      \end{array}
    \right)
    \begin{pmatrix}
      1 & 1-\la & 1-\la \\ 0 & \la & 0 \\ 0 & 0 & \la
    \end{pmatrix}
    \\
    &= \begin{pmatrix}p_n+\la^n & p_n + (1-\la)\la^n & p_n +
      (1-\la)\la^n \\
      q_n & q_n+\la^{n+1} & q_n\\
      r_n & r_n & r_n+\la^{n+1}
    \end{pmatrix}
    \\
    &=
    \left(\begin{array}{lll}p_{n+1}+\la^{n+1} & p_{n+1} & p_{n+1}\\
        q_{n+1} & q_{n+1}+\la^{n+1} & q_{n+1} \\
        r_{n+1} & r_{n+1} & r_{n+1}+\la^{n+1}
      \end{array}\right),
  \end{align*}
  as $p_{n+1}=\sum_0^{n}a_k\la^k=(1-\la)
  \left(\sum_0^{n-1}a_k\la^k+\la^n\right)$, whence
  $p_n+\la^n=p_{n+1}+\la^{n+1}$. For $q_n$ and $r_n$ we have
  $q_{n+1}=q_n,\ r_{n+1}=r_n$. Multiplication by $f_1$ and $f_2$ is
  considered in the same way.
\end{proof}

\begin{cor}\label{cor-image} We have
  \[
  \fe(\De)=
  \begin{cases} x&\ge
    (1-\la)\sum_{k=0}^{n-1}a_k\la^k,\\
    y&\ge
    (1-\la)\sum_{k=0}^{n-1}b_k\la^k,\\
    z&\ge (1-\la)\sum_{k=0}^{n-1}c_k\la^k.
  \end{cases}
  \]
\end{cor}

\begin{proof}The set $\fe(\De)$ is the triangle with the vertices
  $\fe(\p_0), \fe(\p_1)$ and $\fe(\p_2)$. By definition, in this
  triangle $x$ is greater than or equal to the joint first coordinate
  of $\fe(\p_1)$ and $\fe(\p_2)$, i.e., by Lemma~\ref{lem:image1},
  $x\ge (1-\la)\sum_{k=0}^{n-1}a_k\la^k$.  The same argument applies
  to $y$ and $z$.
\end{proof}

\begin{cor}\label{cor-hole} We have
  \[
  \fe(H_0)=
  \begin{cases} x&<
    (1-\la)\left(\la^n+\sum_{k=0}^{n-1}a_k\la^k\right),\\
    y&<
    (1-\la)\left(\la^n+\sum_{k=0}^{n-1}b_k\la^k\right),\\
    z&< (1-\la)\left(\la^n+\sum_{k=0}^{n-1}c_k\la^k\right).
  \end{cases}
  \]
\end{cor}

\begin{proof}The argument is similar to the one in the proof of
  the previous lemma, so we skip it (note that $H_0=\{(x,y,z) : x<
  1-\la, y< 1-\la, z< 1-\la\}$).
\end{proof}

Lemma~\ref{lem:aux1} cannot be reversed in the sense that any
$\fe(H_0)$ is a hole, as we will see in Section~\ref{sec:reverse}.
However, the following assertion shows that once we have one hole,
we  have infinitely many holes.

\begin{lemma}\label{lem:radial} For any $\la\in(1/2,2/3)$ there is
  an infinite number of holes.
\end{lemma}

\begin{figure}
\scalebox{0.7}{\includegraphics{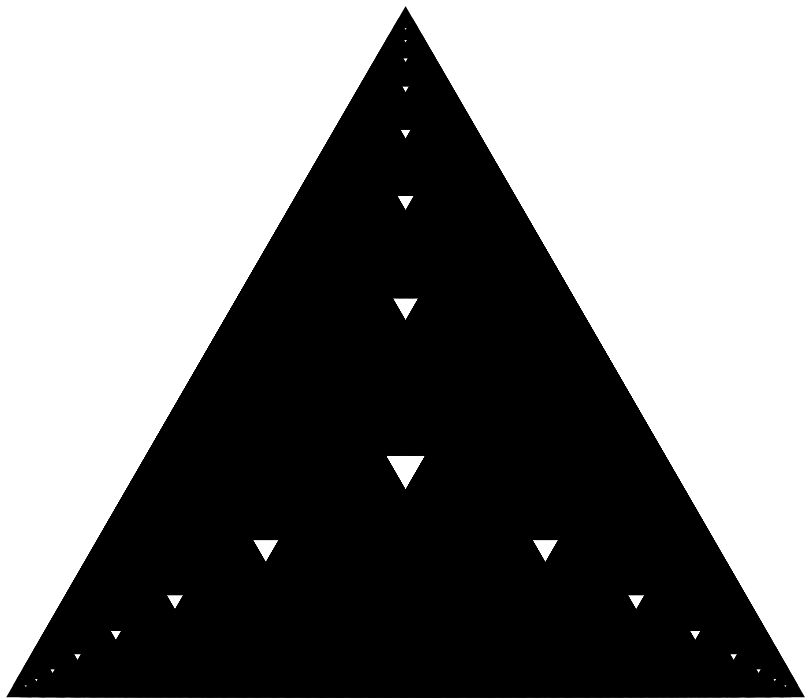}}
\caption{The invariant set for $\la=0.65$.}
\label{fig:radial}
\end{figure}

\begin{proof}We are going to show that $f_i^n(H_0)$ is always a
hole for any $i=0,1,2$ and any $n\ge0$. In view of the symmetry,
it suffices to show that $f_0^n(H_0)$ is a hole. By
Corollary~\ref{cor-hole},
\begin{equation}
f_0^n(H_0)=\{(x,y,z) : x< 1-\la^{n+1},\ y<\la^n(1-\la),\
z<\la^n(1-\la)\}.\label{eq:f0}
\end{equation}
Since the $\De_n$ are nested, it suffices to show that
$f_0^n(H_0)\cap\De_{n+1}=\emptyset$. By Corollary~\ref{cor-image},
this means that the system of inequalities of the form
\begin{equation}
x\ge(1-\la)\sum_0^n a_k\la^k,\ y\ge(1-\la)\sum_0^n b_k\la^k,\
z\ge(1-\la)\sum_0^n c_k\la^k \label{eq:aaa}
\end{equation}
never occurs for $(x,y,z)\in f_0^n(H_0)$. Indeed, if it did, then
by (\ref{eq:f0}), we would have $b_j=c_j=0$ for $0\le j\le n$,
whence $a_0=\dots=a_n=1$, and by (\ref{eq:aaa}),
$x\ge(1-\la)(1+\la+\dots+\la^n)=1-\la^{n+1}$, which contradicts
(\ref{eq:f0}).
\end{proof}
We call any hole of the form $f_i^n(H_0)$ a \emph{radial hole}.
\begin{prop}\label{prop:positive}
Let $\la_*\approx 0.6478$ be the appropriate root of
$$
x^3-x^2+x=\frac12.
$$
Then $\S_\la$ has a nonempty interior if $\la\in[\la_*,2/3)$ and
moreover, each hole is radial---see Figure~\ref{fig:radial}.
\end{prop}

\begin{proof}
We\footnote{We are indebted to B.~Solomyak whose suggestions have
helped us with the idea of this proof.} need to show that for each
$n\ge0$,
\[
F_n:=\De_n\setminus\De_{n+1}\subset\bigcup_{k=0}^n\bigcup_i
f_i^k(H_0)
\]
(within this proof $i$ always runs from 0 to 2). This is obvious
for $n=0,1$, in view of the fact that $H_0$ and the $f_i(H_0)$ are
always holes.

Thus, we have to show that $F_{n}\setminus F_{n-1}$ consists just
of three holes for each $n\ge3$. In view of the symmetry of our
model, this is equivalent to the fact that
\begin{equation*}
f_1 f_0^{n-1}(H_0)\cap F_{n}=\emptyset.
\end{equation*}
We have $F_n=\bigcap_i (\De_n\setminus f_i(\De_n))$. It suffices to
show that $f_1f_0^{n-1}(H_0)\subset f_0(\De_n)$. This in turn follows
from the following relations:
\begin{enumerate}
\item $f_1f_0^{n-1}(H_0)\subset f_0(\De)\cap f_1(\De)$;
\item $f_1f_0^{n-1}(H_0)\cap f_0f_1^{n-1}(H_0)=\emptyset$.
\end{enumerate}
Let $P$ be the vertex of $H_0$ with barycentric coordinates
$(2\la-1,1-\la,1-\la)$. Then (1) is effectively equivalent to
$f_1f_0^{n-1}(P)\in f_0(\De)$, which by Lemma~\ref{lem:image1},
leads to $\la^{n+1}-\la^n+\la\ge\frac12$. By monotonicity of the
root of this polynomial with respect to $n$, the worst case
scenario is $n=2$, which is equivalent to $\la\ge\la_*$.

Let $Q=(1-\la,2\la-1,1-\la)$. The condition~(2) is equivalent to
the fact that the $x$-coordinate of $f_1f_0^{n-1}(P)$ is bigger
than the $x$-coordinate of $f_0f_1^{n-1}(Q)$, which, in view of
Lemma~\ref{lem:image1}, yields the inequality
$\la(1-2\la^{n-1}+2\la^n)>(1-\la)(1+\la^n)$ which is equivalent to
\begin{equation}\label{eq:aux012}
3\la^{n+1}-3\la^n+2\la>1.
\end{equation}
The worst case scenario is $n=3$, where (\ref{eq:aux012}) is
implied by $\la>0.6421$, i.e., well within the range.
\end{proof}

\begin{rmk}
As is easy to see, $\la_*$ is the exact lower bound for the
``purely radial" case, because if $\la<\la_*$, the set
$f_1f_0(H_0)\setminus f_0(\De)$ has an empty intersection with
$\De_3$ and hence is a hole. The details are left to the
interested reader.
\end{rmk}

We finish this section by showing that the boundaries of
$\fe(\De)$ do not contain holes.

\begin{prop}\label{prop:0holes} For $\la\ge 1/2$
$$
\partial{\De}\subset\S_\la.
$$
Consequently, for any $\bm\e$,
$$
\partial{\fe(\De})\subset\S_\la.
$$
\end{prop}

\begin{proof} In  barycentric coordinates,
  $\partial{\De}=\{x=0\}\cup\{y=0\}\cup\{z=0\}$. In view of the
  symmetry, it suffices to show that $K=\{z=0\}\subset\S_\la$. Any
  point of $K$ is of the form $(x,1-x,0)$ with $x\in[0,1]$. Now our
  claim follows from Lemma~\ref{lem:image1} and the fact that every
  $x\in[0,1]$ has the greedy expansion in decreasing powers of $\la$,
  i.e., $x=(1-\la)\sum_1^\infty a_k\la^k$. For $y$ we put $b_k=1-a_k$.

  For the second statement, since $\S_\la$ is invariant,
  $f_i(\S_\la)\subset\S_\la$, whence $\fe(\S_\la)\subset\S_\la$ for
  each $\bm\e$. Now our claim follows from $\partial\fe(\De) =
  \fe(\partial\De)$, together with the first part.
\end{proof}

It follows from this proposition that $\dim_H(\S_\la)\ge1$.

\section{Golden Gaskets}
\label{sec:gg}

Within this section, let $\la$ be equal to the \emph{multinacci
  number} $\om_m$, i.e., the unique positive root of
\[
x^m+x^{m-1}+\dots+x=1,\quad m\ge2.
\]
For every $m$, $\om_m\in(1/2,1)$. In particular, $\om_2$ is the golden
ratio, $\om_2=\frac{\sqrt5-1}2\approx0.618$, $\om_3\approx0.544$, etc.
It is well known that $\om_m\searrow1/2$ as $m\to+\infty$. To simplify
our notation, we simply write $\om$ instead of $\om_m$ within this
section, as our arguments are universal.

We will show that $\S_\om$ is totally self-similar
(Theorem~\ref{thm1}); in Section~\ref{sec:reverse} the converse
will be proved. The key technical assertion is

\begin{prop}The set $\fe(H_0)$ is a hole
for any $\bm\e\in\Si^n$.\label{allholes}
\end{prop}
\begin{proof} Let $\De_n$ be given by (\ref{eq:deltan}), and
\begin{equation}\label{eq:Hn}
H_n:=\bigcup_{\bm\e\in\Si^n}\fe(H_0),\quad n\ge1.
\end{equation}
As in Lemma~\ref{lem:radial}, we  show that
$H_n\cap\De_{n+1}=\emptyset$. By Corollaries~\ref{cor-image} and
\ref{cor-hole}, it suffices to show that the inequalities
\begin{equation}
\begin{aligned}
\om^n + \sum_0^{n-1}a_k\om^k &> \sum_0^{n}\al_k\om^k,\\
\om^n + \sum_0^{n-1}b_k\om^k &> \sum_0^{n}\be_k\om^k,\\
\om^n + \sum_0^{n-1}c_k\om^k &> \sum_0^{n}\ga_k\om^k\label{ineq}
\end{aligned}
\end{equation}
never hold simultaneously, provided all the coefficients are 0's
and 1's, and $a_k+b_k+c_k=\al_k+\be_k+\ga_k=1$.

The key to our argument is the following separation result (we use
the conventional notation here):

\begin{thm}\label{lem:erd} (P. Erd\H os, I. Jo\'o, M. Jo\'o
  \cite[Theorem~4]{EJJ}) Let $\theta>1$, and
  \begin{equation}\label{eq:ltheta}
    \ell(\theta):=
    \inf\,\left\{|\rho|:\rho=\sum_{k=0}^n s_k\theta^k\neq0,
      s_k\in\{0,\pm1\},\ n\ge1\right\}.
  \end{equation}
  Then $\ell(\theta)=\theta^{-1}$ if $\theta^{-1}$ is a multinacci
  number.
\end{thm}

From this lemma we easily deduce a claim about the sums in
question. Indeed, put $\theta=\om^{-1}$ and assume that
$a_k\in\{0,1\},a_k'\in\{0,1\}$ for $k=0,1,\dots,n$, and $\sum_0^n
a_k\om^k > \sum_0^n a'_k\om^k$. Then
\begin{equation}\label{eq:erd}
\sum_0^n (a_k-a_k')\om^k\ge\om^{n+1}
\end{equation}
(just put $s_k=a_{n-k}-a'_{n-k}$).

We use inequality~(\ref{eq:erd}) to improve the
inequalities~(\ref{ineq}). Formally set $a_n=b_n=c_n=1$ and
include the $\om^n$ term of the left hand side of the
inequalities~(\ref{ineq}) with the summation. Then by
(\ref{eq:erd}),
\begin{align*}
\sum_0^{n}a_k\om^k &\ge \sum_0^{n}\al_k\om^k+\om^{n+1},\\
\sum_0^{n}b_k\om^k &\ge \sum_0^{n}\be_k\om^k+\om^{n+1},\\
\sum_0^{n}c_k\om^k &\ge \sum_0^{n}\ga_k\om^k+\om^{n+1},
\end{align*}
which is equivalent to
\begin{equation}
\begin{aligned}
(1-\om)\om^n + \sum_0^{n-1}a_k\om^k &\ge \sum_0^{n}\al_k\om^k,\\
(1-\om)\om^n + \sum_0^{n-1}b_k\om^k &\ge \sum_0^{n}\be_k\om^k,\\
(1-\om)\om^n + \sum_0^{n-1}c_k\om^k &\ge \sum_0^{n}\ga_k\om^k.
\label{ineq2}
\end{aligned}
\end{equation}
By our assumption, just one of the values $\al_n,\be_n,\ga_n$ is
equal to 1. Let it be $\al_n$, say; then the
inequalities~(\ref{ineq2}) may be rewritten as follows:
\begin{align*}
\sum_0^{n-1}a_k\om^k &\ge \sum_0^{n-1}\al_k\om^k+\om^{n+1},\\
(1-\om)\om^n + \sum_0^{n-1}b_k\om^k &\ge \sum_0^{n-1}\be_k\om^k,\\
(1-\om)\om^n + \sum_0^{n-1}c_k\om^k &\ge \sum_0^{n-1}\ga_k\om^k.
\end{align*}
It suffices to again apply (\ref{eq:erd}) to improve the first
inequality. As
$\sum_{k=0}^{n-1}a_k\om^k>\sum_{k=0}^{n-1}\al_k\om^k$, we have
$\sum_{k=0}^{n-1}a_k\om^k-\sum_{k=0}^{n-1}\al_k\om^k\ge\om^n$,
whence
\begin{align*}
\sum_0^{n-1}a_k\om^k &\ge \sum_0^{n-1}\al_k\om^k+\om^n,\\
(1-\om)\om^n + \sum_0^{n-1}b_k\om^k &\ge \sum_0^{n-1}\be_k\om^k,\\
(1-\om)\om^n + \sum_0^{n-1}c_k\om^k &\ge \sum_0^{n-1}\ga_k\om^k.
\end{align*}
Summing up the left and right hand sides, we obtain, in view of
$a_k+b_k+c_k=\al_k+\be_k+\ga_k=1$,
\[
2(1-\om)\om^n+\sum_0^{n-1}\om^k\ge\om^n+\sum_0^{n-1}\om^k,
\]
which implies $\om\le1/2$, a contradiction.
\end{proof}

This claim almost immediately yields the total self-similarity of
the invariant set $\S_\om$:

\begin{thm}\label{thm1} The set $\S_\om$ is totally self-similar in
the sense of Definition~\ref{def:tss}, i.e.,
\[
\fe(\S_\om)=\fe(\De)\cap \S_\om \quad
\mathrm{for}\,\,\mathrm{any}\,\, \bm\e\in\Si^n.
\]
\end{thm}
\begin{proof} Let $H_n$ be defined by (\ref{eq:Hn}).
Since $H_{n+k}=\bigcup_{\bm\e\in\Si^n}\fe(H_k)$, we have
$\fe(H_k)\subset H_{n+k}$. Furthermore, $\fe(H_{k+1})\subset
\fe(\De)$, whence $\fe(H_k)\subset H_{n+k}\cap \fe(\De)$. On the
other hand, by Proposition~~\ref{allholes}, either $\fe(H_0)\cap
f_{\bm\e'}(H_0)=\emptyset$ or $\fe(H_0)=f_{\bm\e'}(H_0)$ for
$\bm\e\in\Si^{n+k}$. Hence the elements of $H_{n+k}$ are disjoint,
and we have
\[
\fe(H_k)=\fe(\De)\cap H_{n+k}.
\]
Since we have proved in Proposition~\ref{allholes} that
$H_{n+k}\cap \De_{n+k-1}=\emptyset$,
\[
\fe(\De_k)=\fe(\De)\cap\De_{n+k-1}.
\]
The claim now follows from the fact that $\De_k\to\S_\om$ in the
Hausdorff metric and from $\fe$ being continuous.
\end{proof}

\section{Dimensions}
\label{sec:dims}

Within this section we continue to assume $\la=\om_m$ for some
$m\ge2$. From Proposition~\ref{allholes} it is easy to show that
$\S_{\om_m}$ is nowhere dense. We  prove more than that:

\begin{prop}\label{prop:leb0}
The two-dimensional Lebesgue measure of $\S_{\om_m}$ is zero.
\end{prop}
\begin{proof} Our proof is based on Theorem~\ref{thm1}. Note
  first that for any measure $\nu$ (finite or not),
\begin{equation*}
\begin{aligned}
\nu(\De)=\nu(&f_0(\De)\cup f_1(\De)\cup f_2(\De)
        \cup H_0)\\
        &-\nu(f_0(\De)\cap f_1(\De))-\nu
        (f_0(\De)\cap f_2(\De))-\nu
        (f_1(\De)\cap f_2(\De)),
\end{aligned}
\end{equation*}
whence by Theorem~\ref{thm1},
\begin{equation}\label{eq-ss}
\begin{aligned}
\nu\S_{\om_m}=\nu(&f_0(\S_{\om_m})\cup f_1(\S_{\om_m})\cup
f_2(\S_{\om_m}))
        \\
        &-\nu (f_0(\S_{\om_m})\cap f_1(\S_{\om_m}))-\nu
        (f_0(\S_{\om_m})\cap f_2(\S_{\om_m}))-\nu
        (f_1(\S_{\om_m})\cap f_2(\S_{\om_m}))
\end{aligned}
\end{equation}
(because $H_0\cap \S_{\om_m}=\emptyset$). The central point of the
proof is that there exists a simple expression for
$f_i(\S_{\om_m})\cap f_j(\S_{\om_m})$ for $i\neq j$. Namely,
\begin{equation}\label{eq:intersect}
f_i(\S_{\om_m})\cap f_j(\S_{\om_m})=f_i f_j^m(\S_{\om_m}).
\end{equation}
To prove this, note first that in view of Theorem~\ref{thm1}, it
suffices to show that
\begin{equation}\label{eq:intt}
f_i(\De)\cap f_j(\De)=f_i f_j^m(\De).
\end{equation}
Moreover, because of the symmetry of our model, in fact, we need
to prove only that $f_0(\De)\cap f_1(\De)=f_0 f_1^m(\De)$. This in
turn follows from Corollary~\ref{cor-image}:
\[
f_0(\De)\cap f_1(\De)=\{(x,y,z) : x\ge 1-\om_m, y\ge 1-\om_m\}
\]
and
\[
f_0 f_1^m(\De)=\{(x,y,z) : x\ge 1-\om_m, y\ge
(1-\om_m)(\om_m+\dots+\om_m^m)=1-\om_m\}.
\]
The relation~(\ref{eq:intersect}) is thus proved. Hence
(\ref{eq-ss}) can be rewritten as follows:
\begin{equation}\label{eq-ss2}
\begin{aligned}
\nu\S_{\om_m}=\nu&(f_0(\S_{\om_m})\cup f_1(\S_{\om_m}))\cup
f_2(\S_{\om_m})
        )\\
        &-\nu (f_0f_1^m(\S_{\om_m}))-\nu
        (f_0f_2^m(\S_{\om_m}))-\nu(f_1f_2^m(\S_{\om_m})).
\end{aligned}
\end{equation}
Finally, let $\nu=\mu$, the two-dimensional Lebesgue measure
scaled in such a way that $\mu(\De)=1$. In view of the $f_i$ being
affine contractions with the same contraction ratio $\om_m$ and by
(\ref{eq-ss2}),
\begin{equation*}
\mu(\S_{\om_m})=3\om_m^2\mu(\S_{\om_m})-
3\om_m^{2(m+1)}\mu(\S_{\om_m}),
\end{equation*}
whence,
\begin{equation}\label{eq:aux1}
(1-3\om_m^2+3\om_m^{2(m+1)})\mu(\S_{\om_m})=0.
\end{equation}
It suffices to show that $1-3\om_m^2+3\om_m^{2(m+1)}\neq0$. For
$m=2$, in view of $\om_2^2=1-\om_2$, this follows from
$1-3\om_2^2+3\om_2^6=\om_2^8>0$; for $m\ge3$, we have
$1-3\om_m^2+3\om_m^{2(m+1)}>1-3\om_m^2>0$, because
$\om_m\le\om_3<0.544<1/\sqrt3$.

Thus, by (\ref{eq:aux1}), $\mu(\S_{\om_m})=0$.
\end{proof}

\begin{rmk}The only fact specific to the multinacci numbers
that we used in this proof is the relation~(\ref{eq:intt}). It is
easy to show that, conversely, this relation implies $\la=\om$. We
leave the details to the reader.
\end{rmk}

We do not know whether the Lebesgue measure of $\S_\la$ is zero if
$\la<\om_2$ (this is what the numerics might suggest), but a
weaker result is almost immediate (NB: $\om_2>1/\sqrt{3}>\om_3$):

\begin{prop}\label{prop:lebmes}
For any $\lambda<1/\sqrt{3}$ the invariant set $\S_\la$ has zero
Lebesgue measure.
\end{prop}

\begin{proof}
Since $\mathcal{S}_\lambda =f_0(\mathcal{S}_\lambda) \cup
f_1(\mathcal{S}_\lambda) \cup f_2(\mathcal{S}_\lambda)$ it follows
that
$$
\mu(\mathcal{S}_\lambda) \le 3\lambda^2 \mu(\mathcal{S}_\lambda).
$$
As $\mathcal{S}_\lambda$ is bounded, we know that
$\mu(\mathcal{S}_\lambda)<\infty$, so that either
$\mu(\mathcal{S}_\lambda)=0$ or $1\le 3\lambda^2$ as required.
\end{proof}

Return to the case $\la=\om_m$. As $\S_{\om_m}$ has zero Lebesgue
measure, it is natural to ask what its Hausdorff dimension is. Let
$\Hh^s$ denote the $s$-dimensional \emph{Hausdorff measure}. As is
well known,
\begin{equation}\label{eq:selfs}
\Hh^s(\la B+\mathbf x)=\la^s \Hh^s(B)
\end{equation}
for any Borel set $B$, any vector $\mathbf x$ and any $\la>0$. Let
us compute $\Hh^s(\S_{\om_m})$. By (\ref{eq-ss2}) and
(\ref{eq:selfs}) with $\nu=\Hh^s$,
\begin{equation*}
\Hh^s(\S_{\om_{m}})=3\om_m^s\Hh^s(\S_{\om_{m}})-
3\om_m^{s(m+1)}\Hh^s(\S_{\om_{m}}).
\end{equation*}
We see that unless
\begin{equation}\label{eq:hd}
1-3\om_m^s+3\om_m^{s(m+1)}=0,
\end{equation}
the $s$-Hausdorff measure of the invariant set is either 0 or
$+\infty$. Recall that the value of $d$ which separates $0$ from
$+\infty$ is called the \emph{Hausdorff dimension} of a Borel set $E$
(notation: $\dim_H(E)$). This argument relies on the invariant set
having non-zero measure in the appropriate dimension, which we do not
know, so in fact only amounts to a heuristic argument suggesting

\begin{thm}\label{thm-hd}
The Hausdorff dimension of the invariant set $\S_{\om_m}$ equals
its box-counting dimension and is given by
$$
\dim_H(\S_{\om_m}) = \dim_B(\S_{\om_m})=
\frac{\log\tau_m}{\log\om_m},
$$
where $\tau_m$ is the largest root of the polynomial
$3z^{m+1}-3z+1$.
\end{thm}

The fact that the Hausdorff dimension coincides with the
box-counting dimension for the attractor of a finite IFS is
universal \cite{Falc}. A rigorous proof of the formula for
$\dim_H(\S_{\om_m})$ is given in the appendix. It amounts to
showing that the invariant set \emph{essentially} coincides with
the invariant set of a countably infinite IFS which satisfies the
OSC.

\begin{rmk}
The case $m=2$ (the golden ratio) is especially nice as here
\[
\tau_2=\frac2{\sqrt3}\cos(7\pi/18).
\]
The authors are grateful to H.~Khudaverdyan for pointing this out.
Note also that there cannot be such a nice formula for $m\ge3$,
because, as is easy to show, the Galois group of the extension
$\mathbb Q(\tau_m)$ with $m\ge3$ is symmetric.
\end{rmk}

\begin{rmk} Let us also mention that the set of holes,
$\De\setminus\S_{\om_2}$, can be identified with the Cayley graph
of the semigroup
\[
\Gamma:=\{0,1,2 \mid 100=011, 200=022, 211=122\},
\]
namely, $f_{\e_0}\dots f_{\e_{n-1}}(H_0)$ is identified with the
equivalence class of the word $\e_0\dots\e_{n-1}$. The relations
$ij^2=ji^2,\,i\neq j$ in $\Gamma$ correspond to the relations
$f_if_j^2=f_jf_i^2,\ i\neq j$.

Thus, $\De\setminus\S_{\om_2}$ may be regarded as a generalization
of the \emph{Fibonacci graph}---the Cayley graph of the semigroup
$\{0,1 \mid 100=011\}$ introduced in \cite{AZ} and studied in
detail in \cite{SV}.

Let $u_n$ stand for the cardinality of level $n$ of $\Gamma$ ($=$
the number of holes of the $n^{\textrm{th}}$ layer). As is easy to
see, $u_0=1, u_1=3, u_2=9$ and
\[
u_{n+3}=3u_{n+2}-3u_n,
\]
whence the rate of growth of $\Gamma$, $\lim_n\sqrt[n]u_n$, is
equal to $\tau_2^{-1}$. This immediately yields another proof that
the box-counting dimension of $\S_{\om_2}$ is equal to its
Hausdorff dimension. The analogous results hold for $\la=\om_m$
for any $m\ge2$. We leave the details to the reader.
\end{rmk}

\begin{table}[ht]
\centerline{
\begin{tabular}{c|c|c}
       $m$ & $\om_m$ & $\dim_H(\S_{\om_m})$ \\
\hline 2 & 0.61803 & 1.93063 \\
\hline 3 & 0.54369 & 1.73219 \\
\hline 4 & 0.51879 & 1.65411 \\
\hline 5 & 0.50866 & 1.61900 \\
\hline 6 & 0.50414 & 1.60201 \\
\hline 7 & 0.50202 & 1.59356 \\
\hline 8 & 0.50099 & 1.58930 \\
\hline 9 & 0.50049 & 1.58715 \\
\hline $\dots$ & $\dots$ & $\dots$ \\
\hline $\infty$ & $1/2$ & $\log3/\log2$
\end{tabular}
}

\vskip0.5truecm

\caption{Hausdorff dimension of $\S_{\om_m}$.}
\label{table1}
\end{table}

\begin{rmk} Recall that $\log3/\log2$ is the Hausdorff
  dimension of the Sierpi\'nski sieve. From Theorem~\ref{thm-hd} it
  follows that $\dim_H(\S_{\om_m})\to\log3/\log2$ as $m\to+\infty$
  (see also Table~\ref{table1}). Thus, although the Hausdorff
  dimension does not have to be continuous, in our case it is
  continuous as $m\to\infty$.
\end{rmk}

\bigskip

\section{The converse and a number-theoretic application}
\label{sec:reverse}

The aim of this section is to show that Theorem~\ref{thm1} can be
reversed, i.e., the choice of multinacci numbers was not accidental.
We are going to need some facts about \emph{$\la$-expansions} of
$x=1$.

Note first that for every $\la\in(1/2,1)$ there always exists a
sequence $(a_k)_1^\infty$ (called a \emph{$\la$-expansion}) that
satisfies
\begin{equation*}
1=\sum_{k=1}^\infty a_k\la^k.
\end{equation*}
The reason why there always some $\la$-expansion available is
because one can always take the \emph{greedy expansion} of 1,
namely, $a_k=[\la^{-1}T_\la^{k-1}(1)]$, where $[\cdot]$ stands for
the integral part, and $T_\la(x)=x/\la-[x/\la]$ (see, e.g.,
\cite{Pa}). We always assume in this preamble that
$\a=(a_k)_1^\infty$ is the greedy $\la$-expansion of 1.

There is a convention in this theory that if the greedy expansion
is of the form $(a_1,\dots,a_N,\linebreak[2] 0,0,\dots)$, then it
is replaced by $(a_1,\dots,a_N-1)^\infty$ (this clearly does not
change the value). For instance, the greedy expansion of 1 for
$\la=\om_2$ is $101010\dots$, and more generally, if $\la=\om_m$,
then $\a=(1^{m-1}0)^\infty$.

\begin{rmk}\label{rmk:1}
As is well known \cite{Pa},
\[
\sum_{k=n+1}^\infty a_k\la^k\le\la^n
\]
for any $n\ge0$, and the equality holds only if $\a$ is purely
periodic, and $a_{n+j}\equiv a_j$ for each $j\ge1$.
\end{rmk}

\begin{lemma}\label{lem:strictlyless}
  Unless $\la$ is a multinacci number, there is always an $n$ such
  that $a_n=0,\ a_{n+1}=1$ and $\sum_{k=n+1}^\infty a_k\la^k<\la^n$.
\end{lemma}

\begin{proof} It follows from Remark~\ref{rmk:1} that unless each
  0 in $\a$ is followed by the string of $L$ 1's for some $L\ge1$
  (which is exactly multinacci), the condition in question is always
  satisfied.
\end{proof}

\begin{thm}\label{thm:converse}
  If, for some $\la\in(1/2,2/3)$, the invariant set $\S_\la$ is
  totally self-similar, then $\la=\om_m$ for some $m\ge2$.
\end{thm}

\begin{proof} Assume $\la$ is such that $\S_\la$ is totally
  self-similar. By definition of total self-similarity, $\fe(H_0)\cap
  S_\la=\emptyset$ for any $\bm\e$, i.e., the claim of
  Proposition~\ref{allholes} must be true. Therefore, it would be
  impossible that, say, $f_0(\De)$ had a ``proper" intersection with
  $\fe(H_0)$ for some $\bm\e$ (see Figure~\ref{fig:triangle})---should
  this occur, a part of $\partial\fe(\De)$ would have a hole, whence
  $\partial\De\not\subset\S_\la$, which contradicts
  Proposition~\ref{prop:0holes}.

  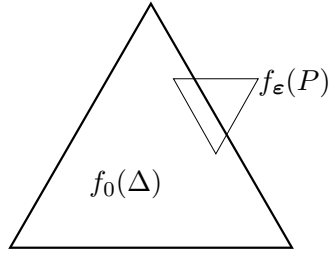
\begin{figure}
    \centering \unitlength=0.75mm
    \begin{picture}(50,45)
      \thicklines
      \path(0,0)(25,43.3)(50,0)(0,0)
      \put(14,10){$f_0(\Delta)$}
      \thinlines
      \path(29,30)(36.5,16.7)(44,30)(29,30)
      \put(44,28){$\fe(P)$}
    \end{picture}
    \caption{The pattern that always occurs unless $\la=\om_m$}
    \label{fig:triangle}
  \end{figure}

  Let us make the necessary computations. Put, as above,
  $P=(2\la-1,1-\la,1-\la)$; then $\fe(P)$ has the $x$-coordinate equal
  to $(2\la-1)\la^n+(1-\la)\sum_0^{n-1}a_k\la^k$ (just apply
  Lemma~\ref{lem:image1}). Assume we have a situation exactly like in
  Figure~\ref{fig:triangle}. As is easy to see,
  $f_0(\De)=\{x\ge1-\la\}$, this $x$-coordinate must be less than
  $1-\la$, whereas the $x$-coordinate of the side that bounds
  $\fe(H_0)$ must be less than $1-\la$. Thus,

  \begin{equation}\label{eq:lessless}
    \frac{2\la-1}{1-\la}\,\la^n < 1-\sum_1^{n-1}a_k\la^k < \la^n
  \end{equation}
  (the sum begins at $k=1$, because obviously $a_k$ must equal 0).
  Thus, we only need to show that if $\la\in(1/2,2/3)$ and not
  multinacci, then there always exists a 0-1 word $(a_1\dots a_{n-1})$
  such that (\ref{eq:lessless}) holds.

  Assume first that $1/2<\la<\om_2$ and not a multinacci number. Let
  $\mathbf a$ be the greedy $\la$-expansion of 1; then by
  Lemma~\ref{lem:strictlyless}, there exists $n\ge1$ such that $a_n=0,
  a_{n+1}=1$, and $1-\sum_0^{n-1}a_k\la^k=\sum_{n+1}^\infty a_k\la^k
  <\la^n$.

  Consider the left hand side inequality in (\ref{eq:lessless}).
  Since $a_{n+1}=1$, we have
  \[
  \sum_{n+1}^\infty a_k\la^k\ge\la^{n+1}>\frac{2\la-1}{1-\la}\,\la^n,
  \]
  as $\la<\om_2$ is equivalent to $\la^2+\la<1$, which implies
  $\la>(2\la-1)/(1-\la)$.

  Assume now $\la>\om_2$ (recall that there are no multinacci numbers
  here). Put $n=2$ and $a_1=1$. Then (\ref{eq:lessless}) turns into
  \[
  \frac{2\la-1}{1-\la}\,\la^2<1-\la<\la^2,
  \]
  which holds for $\la\in(\om_2,\la_*)$, where $\la_*$ is as in
  Proposition~\ref{prop:positive}, i.e., the root of
  $2x^3-2x^2+2x-1=0$.

  Thus, it suffices to consider $\la\in[\la_*,2/3)$. By
  Proposition~\ref{prop:positive}, there are no holes in $f_0(\De)\cap
  f_1(\De)$ at all, which means that $\S_\la$ cannot be totally
  self-similar. The theorem is proved.
\end{proof}

\begin{figure}
\scalebox{0.8} {\includegraphics{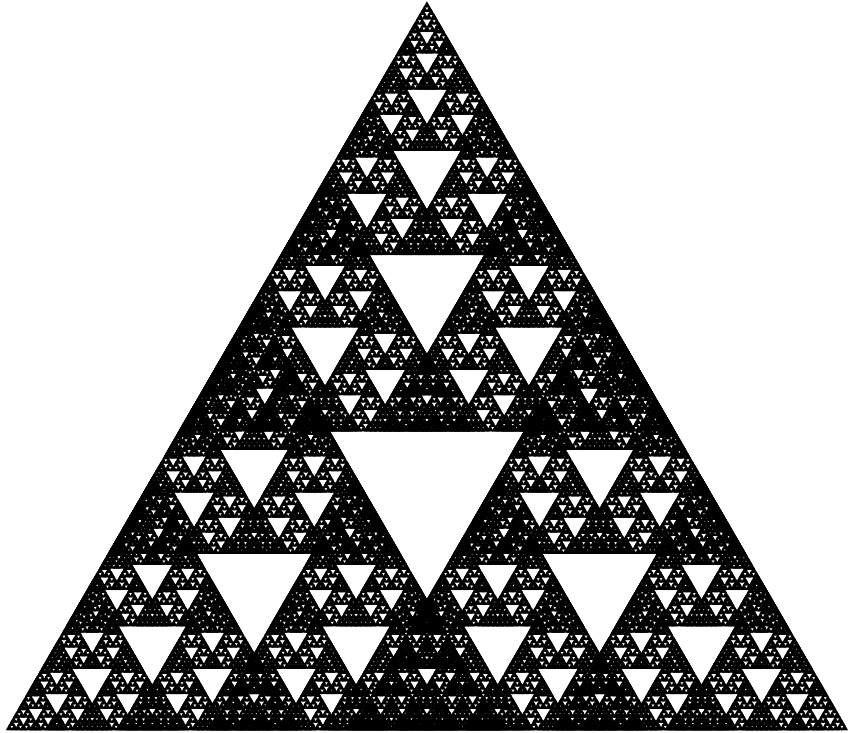}} \caption{The
invariant set for $\la=0.59$. Observe that the holes up to the
second ``layer" seem to be intact, but start to ``deteriorate"
starting with the third ``layer".} \label{fig:059}
\end{figure}

\begin{rmk} Figure~\ref{fig:059} shows consequences of $\S_\la$ being
  not totally self-similar. We see that the whole local structure gets
  destroyed.
\end{rmk}

Theorem~\ref{thm:converse} has a surprising number-theoretic
application (recall the definition of $\ell(\theta)$ is given in
Theorem~\ref{lem:erd}):

\begin{cor}\label{cor:estim}
Let $\theta\in(3/2,2)$. Then either $\theta^{-1}$ is multinacci or
\begin{equation}\label{est}
\ell(\theta)\le\frac{2}{2+\theta}<\frac1{\theta}.
\end{equation}
\end{cor}
\begin{proof} Assume $\la=\theta^{-1}\neq\om_m$ for any $m\ge2$.
From Theorem~\ref{thm:converse} it follows that our method of
proving Proposition~\ref{allholes} simply would not work if $\la$
was not a multinacci number. Recall that our proof was based on
Theorem~\ref{lem:erd} which must consequently be wrong if $\la$ is
not multinacci.

Moreover, with $\kappa:=\theta \ell(\theta)$ and by the same chain of
arguments as in the proof of Proposition~\ref{allholes}, we come
at the end to the inequality
\[
2(1-\kappa\la)\ge\kappa
\]
(in the original proof we had it with $\kappa=1$) which is equivalent
to $\ell(\theta)\le\frac{2}{2+\theta}$. Thus, if this inequality is
\emph{not} satisfied, then the system of inequalities~(\ref{ineq})
does not hold either, which leads to the conclusion of
Proposition~\ref{allholes} and consequently yields
Theorem~\ref{thm1}---a contradiction with Theorem~\ref{thm:converse}.
\end{proof}

\begin{rmk} As is well known since the pioneering work \cite{Ga},
  if $\theta$ is a \emph{Pisot number} (an algebraic integer $>1$
  whose Galois conjugates are all less than 1 in modulus), then
  $\ell(\theta)>0$ (note that the $\om_m^{-1}$ are known to be Pisot).
  Furthermore, if $\theta$ is not an algebraic number satisfying an
  algebraic equation with coefficients $0,\pm1$, then by the
  pigeonhole principle, $\ell(\theta)=0$. There is a famous conjecture
  that this is also true for all algebraic non-Pisot numbers.

  Thus (modulo this conjecture), effectively, the result of
  Corollary~\ref{cor:estim} is of interest if and only if $\theta$ is
  a Pisot number. The restriction $\theta>3/2$ then is not really
  important, because in fact, there are only four Pisot numbers below
  $3/2$, namely, the appropriate roots of $x^3=x+1$ (the smallest
  Pisot number), $x^4=x^3+1$, $x^5-x^4-x^3+x^2=1$ and
  $x^3=x^2+1$.\footnote{In fact, there is just a finite number of
  Pisot numbers below $\frac{1+\sqrt5}2$, and they all are known
  \cite{Pisot}.} The respective values of $\ell(\theta)$ for these
  four numbers are approximately as follows: $0.06, 0.009, 0.002,
  0.15$ (see \cite{BH}), i.e., significantly less than the
  estimate~(\ref{est}).
\end{rmk}

Thus, we have proved

\begin{prop} For each Pisot number $\theta\in(1,2)$ that does not
satisfy $x^m=x^{m-1}+x^{m-1}+\dots+x+1$ for some $m\ge2$,
\[
\ell(\theta)\le\frac{2}{2+\theta}.
\]
\end{prop}

For the history of the problem and the tables of $\ell(\theta)$ for
some Pisot numbers $\theta$ see \cite{BH}.

\begin{rmk} We are grateful to K.~Hare who has indicated the
paper \cite{Z} in which it is shown that $l(q)<2/5$ for
$q\in(1,2)$ and $q^{-1}$ not multinacci. This is stronger than
(\ref{est}) but the proof in \cite{Z} is completely different,
rather long and technical, so we think our result is worth
mentioning.
\end{rmk}

\section{Higher-dimensional analogues}
\label{sec:higher}

The family of IFSs we have been considering consists of three
contractions in the plane, with respective fixed points at the
vertices of a regular 3-simplex in $\BR^2$.  In $\BR^d$ it is
natural to consider $d+1$ linear contractions with fixed points at
the vertices of the $d+1$-simplex:
$$
f_j(\mathbf x) = \lambda\mathbf x +
(1-\lambda)\p_j,\quad(j=0,\dots,d).
$$
For example, when $d=3$ the four maps are contractions towards the
vertices of a regular tetrahedron in $\BR^3$.

Using the analogous barycentric coordinate system ($x_j$ is the
distance to the $j^{\textrm{th}}$ $(d-1)$-dimensional face of the
simplex), the maps $f_0,\dots, f_d$ are given by matrices
analogous to those in Lemma~\ref{lem:matr}. The algebra of these
maps is directly analogous to the family of three maps we have
considered so far. The proofs of the following results are left as
exercises (most are extensions of corresponding results earlier in
the paper).

\medskip

\noindent \textbf{(1) } If $\la\in[\frac{d}{d+1},1)$, then
$\S_\la=\De$, so there are no holes in the attractor.

\smallskip

\noindent \textbf{(2) } If $\la\le1/2$ the IFS satisfies the Open
Set Condition, and the invariant set is self-similar with
Hausdorff dimension
$$
\dim_H(\mathcal{S}_\lambda) = \frac{\log(d+1)}{-\log\lambda}.
$$

\smallskip

\noindent \textbf{(3) } Since the $(d+1)$-simplex contains the
$d$-simplex at each of its faces, for any fixed $\lambda$ we have
$\S_\lambda(d+1)\supset\S_\lambda(d)$ and consequently,
$$
\dim_H(\mathcal{S}_\lambda(d+1)) \ge
\dim_H(\mathcal{S}_\lambda(d)).
$$

\smallskip

\noindent \textbf{(4) }
If $\lambda=\omega_m$ (the multinacci number), then the invariant set
is totally self-similar, and the dimension $s$ satisfies
$$
s = \frac{\log\tau_{m,d}}{\log\om_m}
$$
where $\tau_{m,d}$ is the largest root of $\frac12
d(d+1)t^{m+1}-(d+1)t +1 =0$.  See Table~\ref{tab:dimensions} for some
values.  One can see from this that, for fixed $m$ and large $d$, the
Hausdorff dimension increases logarithmically in $d$.

\smallskip

\noindent \textbf{(5) }
If $\lambda<(d+1)^{-1/d}$, then---similarly to
Proposition~\ref{prop:lebmes}---$\S_\la$ has zero $d$-dimension\-al
Lebesgue measure, but we do not know what happens for $\la\in
\bigl((d+1)^{-1/d},\frac{d}{d+1}\bigr)$.

\begin{table}[ht]
  \centering
  \begin{tabular}{l|l|l|l|l|l|c|l}
    $d$ & $\omega_2$ & $\omega_3$ & $\omega_4$ & $\omega_5$ &
    $\omega_6$ & \dots & 1/2\\
    \hline
    2 &  1.93 & 1.73 & 1.65 & 1.62 & 1.60 & \dots& 1.583\\
   \hline 3 &  2.61 & 2.23 & 2.10 & 2.05 & 2.02 & \dots& 1.999 \\
   \hline 4 &  3.13 & 2.61 & 2.45 & 2.38 & 2.35 & \dots& 2.322 \\
   \hline 5 &  3.54 & 2.92 & 2.72 & 2.65 & 2.62 & \dots& 2.585 \\
   \hline 6 &  3.89 & 3.18 & 2.96 & 2.88 & 2.84 & \dots& 2.807
  \end{tabular}
  \vskip0.5truecm
  \caption{Hausdorff dimension of golden $d$-gaskets}
  \label{tab:dimensions}
\end{table}

\section{Final remarks and open questions}
\label{sec:final}

\noindent {\bf (1)} The fact that the triangle is equilateral in
our model is unimportant. Indeed, given any three non-collinear
points $\p'_0,\p'_1,\p'_2$ in the plane there is a (unique) affine
map $A$ that maps each $\p'_j$ to the corresponding $\p_j$ we have
been using. For given $\la$ let $\S'_\la$ be the invariant set of
the IFS defined by (\ref{eq:ifs0}) with the $\p'_j$ in place of
the $\p_j$.  Then it is clear that $\S_\la = A(\S'_\la)$.  For a
given value of $\la$ all the invariant sets are therefore affinely
equivalent, and in particular have the same Hausdorff dimension
(when this is defined).

\smallskip
\noindent {\bf (2)} The sequence of golden gaskets $\S_{\om_m}$
provides confirmation of some observations regarding the dimension
of fractal sets generated by IFSs where the Open Set Condition
fails.  In particular, a theorem of Falconer \cite{falc88} states
that given linear maps $T_1,\dots,T_k$ on $\mathbb{R}^n$ of norm
less than $1/3$, there is a number $\delta$ such that the
invariant set $F(a_1,\dots,a_k)$ of the IFS $\{T_1+a_1,\dots,
T_k+a_k\}$ has Hausdorff dimension $\delta$ for a.e. $(a_1,\dots,
a_k)\in \mathbb{R}^{nk}$.  In the case that the $T_j$ are all the
same similarity by a factor of $\la$, the dimension is given by
$\delta = \delta(\la) = -\log k/\log\la$.

It has been pointed out \cite{Sol} that the upper bound $1/3$ can
be replaced by $1/2$, but that the theorem fails if the upper
bound is replaced by $1/2+\e$ for any $\e>0$.  This can also be
seen from the golden gaskets $\S_{\om_m}$: given $\e>0$ there is
an $m$ such that $1/2<\om_m<1/2+\e$, and the dimension of the
invariant set $\dim_H(\S_{\om_m})<\delta(\om_m)$.

\smallskip
\noindent {\bf (3)} If one endows each of the maps $f_i$ with
probability $1/3$, this yields a \emph{probabilistic IFS}. Its
general definition can be found, for example, in the survey
\cite{DF}. Then $\S_\la$ becomes the support for the
\emph{invariant measure}; the question is, what can be said about
its Hausdorff dimension? In particular, we conjecture that,
similarly to the 1D case (see \cite{AZ, SV}), it is strictly less
than $\dim_H(\S_\la)$ for $\la=\om_m$.

\smallskip
\noindent {\bf (4)} The main problem remaining is to determine for
which $\la$ the attractor $\S_\la$ has positive Lebesgue measure
and for which zero Lebesgue measure. The numerics suggests the
following

\medskip\noindent\textbf{Conjecture.} (1) For each $\la\in(\om_2,2/3)$
the attractor $\S_\la$ has a nonempty interior (recall that we know
this for $\la\in[\la_*,2/3)$---Proposition~\ref{prop:positive}).
\newline (2) For each $\la\in(1/\sqrt3,\om_2)$ it has an empty
interior and possibly zero Lebesgue measure.

\smallskip
\noindent {\bf (5)} The same range of problems can be considered
for any collection of similitudes $f_j(\mathbf x)=\la\mathbf x
+(1-\la)\mathbf p_j$ in $\BR^d$, where the $\mathbf p_j$ are
vertices of a (convex) polytope~$\Pi$. For instance, are there any
totally-self similar attractors if $\Pi$ is not a simplex and the
OSC fails? This question seems to be particularly interesting if
$d=2$ and $\Pi$ is regular $n$-gon with $n\ge5$.

\section*{Appendix}
\label{sec:appendix}

\addtocounter{section}{1} \setcounter{section}{1}
\def\thesection{\Alph{section}}

We now give a rigorous proof of Theorem~\ref{thm-hd} (repeated
below for convenience), using the fact that the invariant set
almost coincides with the invariant set for an infinite IFS which
satisfies the open set condition, and relying on some results
about such systems \cite{MU}.  We begin with an elementary lemma.
Recall that the multinacci number $\om_m$ is the unique root of
$t^{m+1}-2t+1$ lying in $(\frac12,\frac23)$.

\begin{lemma}\label{lem-auxapp}
  For each integer $m\ge2$, let $\tau_m\in(0,1/2)$ be the smaller
  positive root of $3t^{m+1}-3t+1$, and $\si_m\in(0,1)$ the smaller
  positive root of $2t^m-3t+1$. Then
  $$\textstyle\frac13 < \tau_m < \si_m < \om_m < \frac23.$$
  Consequently
  $$\frac{\log\tau_m}{\log\om_m} > \frac{\log\si_m}{\log\om_m} > 1.$$
\end{lemma}

\begin{proof}
  Let $p_m = 3t^{m+1}-3t+1$ and $q_m =2t^m-3t+1$.  Notice that the
  derivatives of $p_m$ and $q_m$ are monotonic on the interval
  $[0,1]$, so that each have at most two roots on that interval. Note
  also that $p_m(\om_m)<0$ and $q_m(\om_m)<0$. Since $p_m(1)=1$ and
  $q_m(1)=0$ and $p_m(\frac13)>0$ and $q_m(\frac13)>0$ it follows that
  $\frac13<\tau_m<\om_m$ and $\frac13<\si_m<\om_m$.

  Finally, $p_m(\si_m) = 3\si_m \left( \frac{3\si_m-1}{2} \right) -
  3\si_m + 1 = \frac12(3\si_m-1)(3\si_m-2)<0$, so that $\si_m>\tau_m.$
\end{proof}

To complete the picture (though we don't use this), if $\tau_m'$
and $\si_m'$ are the other positive roots of $p_m$ and $q_m$
respectively, then $\om_m<\tau_m'<\si_m'=1$. Furthermore,
$\lim_{m\to\infty}\tau_m = \lim_{m\to\infty}\si_m =\frac13$.

\setcounter{section}{4} \setcounter{lemma}{3}
\def\thesection{\arabic{section}}

\begin{thm}
The Hausdorff dimension of the invariant set $\S_{\om_m}$ is given
by
$$
\dim_H(\S_{\om_m}) = \frac{\log\tau_m}{\log\om_m},
$$
where $\tau_m$ is defined in the lemma above.
\end{thm}

\setcounter{section}{1} \setcounter{lemma}{1}
\def\thesection{\Alph{section}}

\begin{Def}An alternative definition of $\S_\la$ is as follows
  (see, e.g., \cite{DF}): to any $\bm\e\in\Si^\infty$ there
  corresponds the unique point $\xe = \lim_{n\to\infty}f_{\e_0}\dots
  f_{\e_n}(\mathbf x_0)\in \S_\la$. This limit does not depend on the
  choice of $\mathbf x_0$; we call $\bm\e$ an \emph{address} of $\xe$.
  Note that a given $\mathbf x\in\S_\la$ may have more than one
  address---see Proposition~\ref{prop:contin} below.
\end{Def}

\begin{Def}Let $\U_\la$ denote the \emph{set of uniqueness}, i.e.,
  \[
  \U_\la=\left\{\mathbf x\in\De \mid \exists ! (\e_0,\e_1,\dots) :
    \mathbf x = \xe\right\}.
  \]
In other words, $\U_\la$ is the set of points in $\S_\la$, each of
which has a unique address. These sets seem to have an interesting
structure for general $\la$'s, and we plan to study them in
subsequent papers. Note that in the one-dimensional case
($f_j(x)=\la x+(1-\la)j,\ j=0,1$) such sets have been studied in
detail by P.~Glendinning and the third author in \cite{GS}.
\end{Def}

In the course of the proof of this theorem, we also prove the
following

\begin{thm}\label{thm:dimU}
The set of uniqueness $\U_{\om_m}$ is a self-similar set of
Hausdorff dimension
$$
\dim_H(\U_{\om_m}) = \frac{\log\si_m}{\log\om_m},
$$
where $\si_m$ is defined in Lemma~\ref{lem-auxapp}. In particular,
$\si_2=1/2$.
\end{thm}

\begin{proof}[Proof of both theorems]
  The proof proceeds by showing that there is another IFS (an infinite
  one) which does satisfy the OSC, and whose invariant set
  $\A_{\om_m}$ satisfies $\S_{\om_m}=\A_{\om_m}\cup\,\U_{\om_m}$, with
  $\dim_H(\U_{\om_m}) < \dim_H(\A_{\om_m})$.  It then follows that
  $\dim_H(\S_{\om_m}) = \dim_H(\A_{\om_m})$, and the latter is given
  by a simple formula.

  The proof for $m=2$ differs in the details from that for $m>2$ so we
  treat the cases separately.  Note that within this appendix we
  assume the triangle $\De$ has unit side.

  \medskip

  \noindent\textbf{The case $\bm{m}\mathbf{=2}$. }
  Refer to Figure~\ref{fig:trapezia om2} for the geometry of this
  case. Begin by removing from the equilateral triangle $\De$ the
  (open) central hole $H_0$, the three (closed) triangles of side
  $\om_2^2$ that are the images of the three $f_j^2$ ($j=0,1,2$), and
  three smaller triangles of side $\om_2^3$ that are the images of
  $f_if_j^2=f_jf_i^2$.  This leaves three trapezia, whose union we
  denote $T_1$. See Figure~\ref{fig:trapezia om2}~(a).  For this part
  of the proof we write $F_k=f_if_j^2$ (where $i,j,k$ are distinct).

  \begin{figure}[t]
    \centering \unitlength=1.3mm
    \begin{picture}(50,53)(0,-7)
      \thinlines
      \dottedline(0,0)(25,43.3)(50,0)(0,0)
      \thicklines

\path(19,0)(34.55,26.85)(15.45,26.85)(31,0)(40.5,16.45)(9.5,16.45)(19,0)
      \put(31,34){$\omega_2^2$}
      \put(38,22){$\omega_2^3$}
      \put(25,-8){\hbox to 0pt{\hss(a) 3 trapezia\hss}}
    \end{picture}
    \qquad
    \begin{picture}(50,53)(30,13)  
      \thicklines
      \path(30.902, 53.523)(69.098, 53.523)(57.295,
33.079)(42.705,33.079)(30.902, 53.523)
      \path(54.509,53.523)(42.705, 33.079)
      \path(45.492,53.523)(57.295, 33.079)
      \put(48,29.5){$\omega_2^n$}
      \put(52,48){$\omega_2^{n+1}$}
      \put(30,17){
          \begin{minipage}[m]{0.35\textwidth}
            (b) decomposition of a trapezium at level $n$: 1 hole, 1
            triangle and 2 (smaller) trapezia
          \end{minipage}}
    \end{picture}

    \caption{Decomposing the golden gasket ($\lambda=\omega_2$)}
    \label{fig:trapezia om2}
  \end{figure}
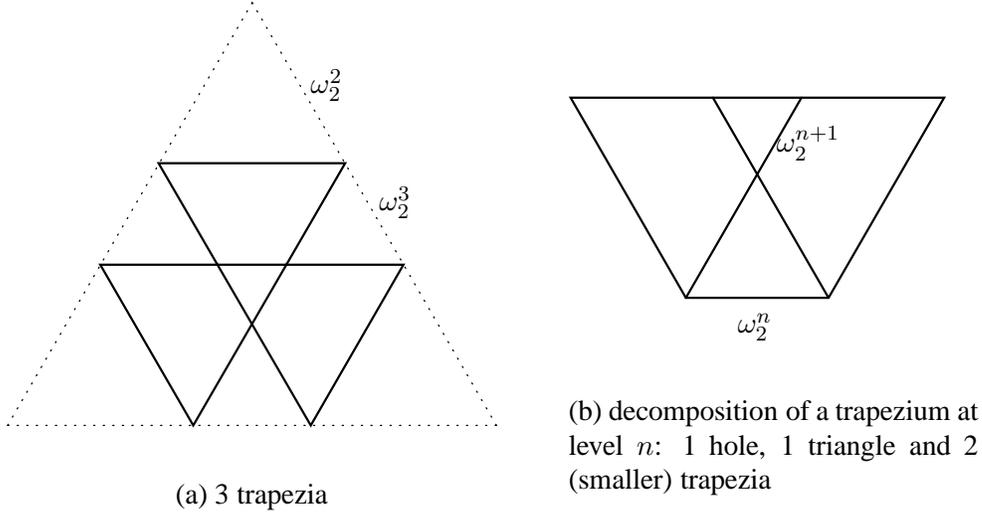

  Each of the three trapezia is decomposed into the following sets: a
  hole (together forming $H_1$), an equilateral triangle of side
  $\om_2^4$, and two smaller trapezia---smaller by a factor of
  $\om_2$.  The three equilateral triangles at this level are the
  images of $f_1f_0f_2^2=f_1F_1$ (for the lower left trapezium),
  $f_2F_2$ (lower right) and $f_0F_0$ (upper trapezium). See
  Figure~\ref{fig:trapezia om2}~(b).  At the next level the
  equilateral triangles are the images of $f_jf_iF_i$ with $i\neq j$,
  and at the following $f_kf_jf_iF_i$ with $k\neq j$ and $j\neq i$.

  This decomposition of the trapezia is now continued \emph{ad
  infinitum}. At the $n^{\textrm{th}}$ level there are
  $3\cdot2^{n-1}$ holes forming $H_n$, the same number of equilateral
  triangles that are images of similarities by $\om_2^n$ and twice as
  many trapezia. Note that at each stage, the holes consist of those
  points with no preimage, the equilateral triangles of those points
  with two preimages and the trapezia of points with a unique
  preimage.

  Let $\A_{\om_2}$ be defined as the attractor corresponding to the
  equilateral triangles in the above construction; thus, it is the
  attractor for the infinite IFS with generators
  \begin{equation}
    \label{eq:IFS2}
    \left\{f_j^2, F_j, f_jF_j, f_jf_iF_i, f_kf_jf_iF_i,\dots\right\},
  \end{equation}
  where the general term is of the form $f_{j_1}f_{j_2}\dots
  f_{j_n}F_{j_n}$ with adjacent $j_k$ different from each other.
  Notice that this IFS satisfies the open set condition. In \cite{MU}
  a deep theory of \emph{conformal IFS} (which our linear one
  certainly is) has been developed. From this theory it follows that,
  similarly to the finite IFSs, the Hausdorff dimension $s$ of the
  invariant set $\A$ (henceforward we drop the subscript `$\om_2$')
  equals its {\rm similarity dimension} given by $\om_2^s=\tau$, where
  in our case, $\tau$ is a solution of
  $$
  1=3\tau^2 + 3\tau^3 + 3\tau^4\sum_{n=0}^\infty 2^n\tau^{n}.
  $$
  This equation has a unique positive solution and is equivalent to
  $(3\tau^3-3\tau+1)(\tau+1)=0$ provided $\tau<1/2$ (the radius of
  convergence of the above power series). Thus, $\tau$ is the solution
  of
  \begin{equation*}
    3\tau^3-3\tau+1=0,
  \end{equation*}
  with $\tau<1/2$, in agreement with the value of the dimension of
  $\S$ given in the theorem.

  Since this IFS is contained in the original IFS (generated by the
  $f_i$), so $\A\subset \S$.

  \begin{figure}
    \centering \scalebox{0.8} {\includegraphics{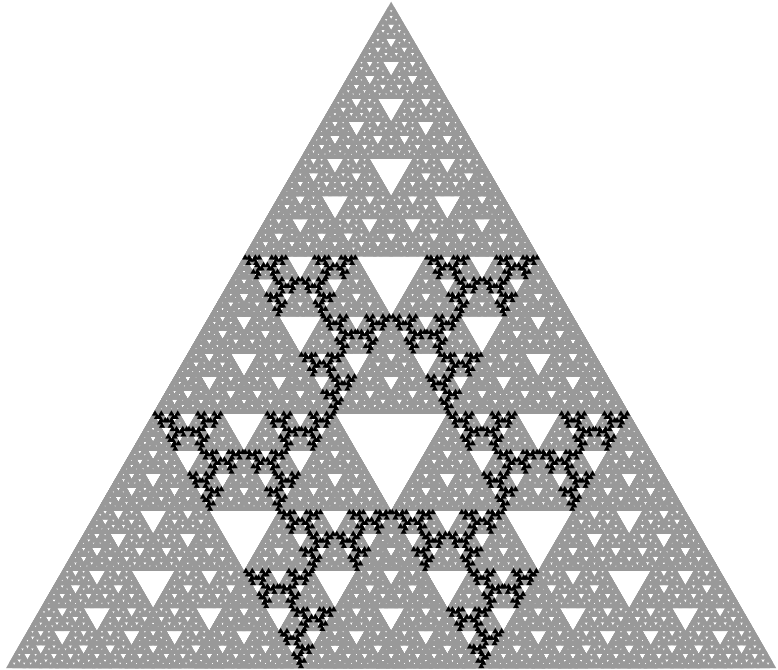}}
    \caption{The fractal set $\U'_{\om_2}$ superimposed on a grey
      $\S_{\om_2}$.}
    \label{fig:u2}
  \end{figure}

  Now let $\U'$ be the limit of the sequence of unions of trapezia
  defined by the above procedure: write $\U^{(n)}$ for the union of
  the $3\cdot2^{n-1}$ trapezia obtained at the $n^{\textrm{th}}$ step,
  then $\U^{(n+1)}\subset \U^{(n)}$ and
  $$
  \U' = \bigcap_{n>0} \U^{(n)}.
  $$
  By construction, $\U'$ is a connected self-similar Cantor set,
  with dimension
  \begin{equation}
    \label{eq:dimU2}
    \dim_H(\U')=-\frac{\log2}{\log\om_2}.
  \end{equation}
  This follows from the standard arguments, since $\#\, \U^{(n)}
  \asymp 2^n$ and $\mbox{diam}\, \U^{(n)} \asymp \om_2^n$.

  We claim that
  $$
  \U'\cup \bigcup_{n=1}^\infty\bigcup_{k=0}^2 f_k^{2n}(\U')=\U.
  $$
  To see this, we turn to the addresses in the symbol space $\Si$.
  In view of the relation $f_if_j^2=f_jf_i^2$, each $\mathbf x\in\S$
  that has multiple addresses, must have
  $\e_{j-1}\neq\e_j,\e_j=\e_{j+1}$ for some $j\ge1$. By our
  construction, this union is the set of $\mathbf x$'s whose addresses
  can have equal symbols only at the beginning.  Thus, it is indeed
  the set of uniqueness.

  Since this is a countable union of sets of the same dimension, it
  follows that $\dim_H\U = \dim_H\U'$ (see, e.g., \cite{Falc}). We
  claim that
  \begin{equation}
    \label{eq:S=AU}
    \S = \A \cup\U
  \end{equation}
  with the union being disjoint.  The result (for $\om_2$) then
  follows, since $\dim\A>\dim\U$.

  To verify (\ref{eq:S=AU}), we observe that from (\ref{eq:IFS2}) it
  follows that $\mathbf x\in\A$ if and only if $\mathbf x = \xe$ for
  some $\bm\e$ for which there are two consecutive indices which
  coincide and do not coincide with the previous one, i.e.,
  $\A\cap\U=\emptyset$. Conversely, every point in $\S$ with more than
  one address lies in $\A$. Thus, apart from the three vertices of
  $\De$, $\S\setminus \A$ consists of the points with a unique
  address, and expression (\ref{eq:dimU2}) proves
  Theorem~\ref{thm:dimU} for the case $m=2$.

  \medskip

  \noindent\textbf{The case $\bm{m}\mathbf{\ge3}$. }
  The overall argument is similar to that for $m=2$, except that the
  trapezia are replaced by hexagons, and the recurrent structure is
  consequently different (more complicated); we only describe where
  the arguments differ.

  We begin in the same way, by removing the central hole $H_0$, and
  decomposing the remainder into 3 small triangles of side $\om_m^m$
  at the vertices---the images of $f_j^m$, 3 smaller triangles of side
  $\om_m^{m+1}$ on each side---images of $f_i f_j^m$, and 3 remaining
  hexagons (instead of trapezia).

  These hexagons have sides of length $\om_m^m,\, (1-\om_m-\om_m^m),
  \,\om_m^{m+1}, \, (2-3\om_m), \, \om_m^{m+1}$ and
  $(1-\om_m-\om_m^m)$ (in cyclic order). We call hexagons similar to
  these, $\om$-hexagons, and this one in particular an $\om$-hexagon
  of \emph{size} $\om^m$.  Notice that these $\om$-hexagons have a
  single line of symmetry, and the \emph{size} refers to the length of
  the smaller of the two sides that meet this line of symmetry.

  Each $\om$-hexagon of size $\om_m^k$ can be decomposed into: $(m-1)$
  holes of various sizes down the line of symmetry; $(3m-5)$
  equilateral triangles, 3 each of sizes $\om_m^{k+2},
  \om_m^{k+3},\dots,\om_m^{k+m-1}$ and one of size $\om_m^{k+m}$; this
  leaves $2(m-1)$ $\om$-hexagons, 2 each of sizes $\om_m^{k+1},
  \om_m^{k+2}, \dots, \om_m^{k+m-1}$ (see Figure \ref{fig:hexagon
    decomposition} for the cases $m=3$ and $4$).

  \begin{figure}[t]
    \centering \unitlength=1.4mm
    \begin{picture}(50,53)(30,30)  
      \thicklines
      \path(41.965, 72.684)(58.036, 72.684)(72.816, 47.085)(68.447,
      39.517)(31.554, 39.517)(27.185, 47.085)(41.965, 72.684)
      \path(37.596, 65.117)(52.376, 39.517)
      \path(62.404, 65.117)(47.624, 39.517)
      \path(35.220, 61.002)(64.780, 61.002)
      \path(44.548, 72.684)(51.292, 61.003)
      \path(55.453, 72.684)(48.709, 61.003)
      \put(49,75){$\omega_3^k$}
      \put(30,43){$\omega_3^{k+1}$}
      \put(71,42){$\omega_3^{k+1}$}
      \put(50,35){\hbox to 0pt{\hss(a)\quad $m=3$\hss}}
      \put(50,31){\hbox to 0pt{\hss 2 holes, 4 triangles and 4
          $\om$-hexagons\hss}}
    \end{picture}
    \quad
    \begin{picture}(50,53)(25,32) 
      \thicklines
      \path(46.378, 80.329)(53.622, 80.329)(74.061, 44.928)(72.182,
      41.674)(27.819, 41.674)(25.940, 44.928)(46.378, 80.329)
      \path(37.518, 64.982)(50.975, 41.674)
      \path(62.483, 64.982)(49.026, 41.674)
      \path(36.543, 63.294)(63.457, 63.294)
      \path(43.524, 75.385)(50.506, 63.293)
      \path(56.476, 75.385)(49.495, 63.293)
      \path(43.019, 74.510)(56.982, 74.510)
      \path(46.903, 80.329)(50.263, 74.510)
      \path(53.097, 80.329)(49.738, 74.510)
      \put(49,82){$\omega_4^k$}
      \put(20,41){$\omega_4^{k+1}$}
      \put(74,42){$\omega_4^{k+1}$}
      \put(50,37){\hbox to 0pt{\hss(b)\quad $m=4$\hss}}
      \put(50,33){\hbox to 0pt{\hss 3 holes, 7 triangles and 6
          $\om$-hexagons\hss}}
    \end{picture}
    \caption{Decomposing an $\om$-hexagon of size $\om_m^k$}
    \label{fig:hexagon decomposition}
  \end{figure}
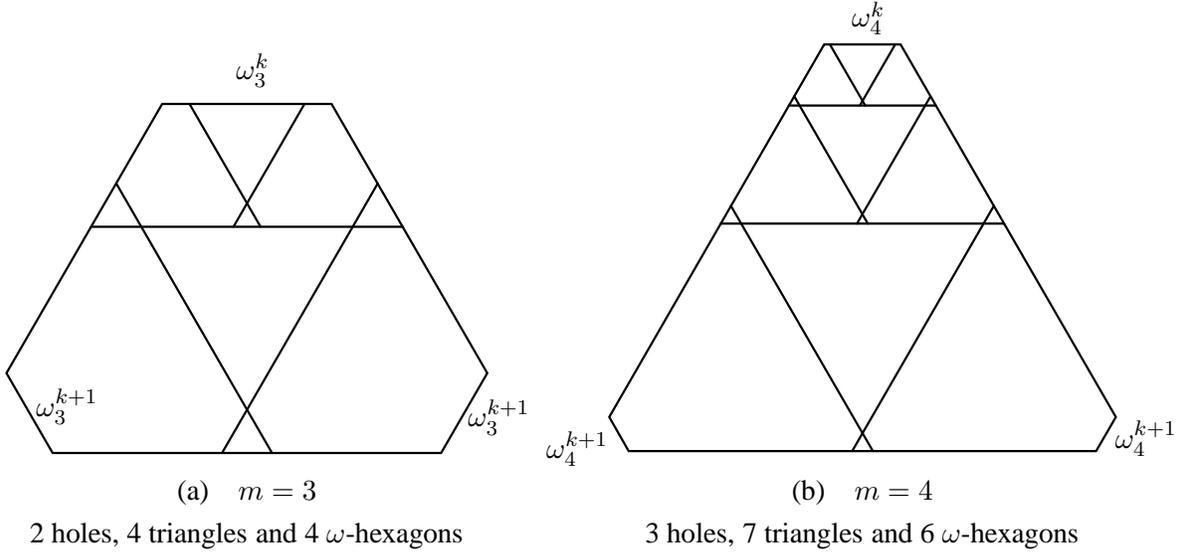

  In the same way as in the case $m=2$, the equilateral triangles
  occurring in this decomposition are the images of the original
  triangle $\De$ under certain similarities arising in the IFS
  generated by $\{f_0,f_1,f_2\}$. This sub-IFS defines a countable IFS
  which satisfies the OSC, permitting us again to compute the
  dimension of the corresponding invariant set $\A_{\om_m}$. We use
  generating functions to compute this dimension.

  Each hexagon of size $\om_m^k$ decomposes into 2 hexagons of sizes
  $\om_m^{k+1},\dots,\om_m^{k+m-1}$.  Thus, each hexagon of size
  $\om_m^k$ arises from decomposing hexagons of sizes
  $\om_m^{k-m+1}\dots \om_m^{k-1}$.  Let $h_k$ be the number of
  hexagons of size $\om_m^k$ that appear in the procedure. Then,
  $h_k=0$ for $k<m$, $h_m=3$ and for $k>m$,
  $$
  h_k = 2\left(h_{k-m+1} + \cdots + h_{k-1}\right).
  $$
  Applying the usual generating function approach, let $Q =
  \sum_{k=1}^\infty h_kt^k$. Then
  \begin{eqnarray*}
    Q &=& 3t^m + 2\sum_{k=m+1}^\infty \sum_{r=1}^{m-1}h_{k-r}t^k \\
    &=& 3t^m + 2\sum_{r=1}^{m-1}t^r\sum_{k=m+1}^\infty h_{k-r}t^{k-r} \\
    &=& 3t^m +  2Q\sum_{r=1}^{m-1}t^r.
  \end{eqnarray*}
  Finally, provided $|t|<r_m$ the radius of convergence of the power
  series,
  $$
  Q = \frac{3t^m(1-t)}{1-3t+2t^m}.
  $$
  Note from its definition in Lemma \ref{lem-auxapp} that
  $r_m=\si_m$.  Now for the triangles: each $\om$-hexagon of size
  $\om_m^k$ gives rise to 3 triangles of sizes
  $\om_m^{k+2},\dots,\linebreak[1]\om_m^{k+m-1}$ and one of size
  $\om_m^{k+m}$. Let there be $p_k$ triangles of size $\om_m^k$.  Then
  $p_k=0$ for $k<m$, $p_m=p_{m+1}=3$, and for $k>m+1$,
  $$
  p_k = h_{k-m}+3(h_{k-m+1}+\cdots+ h_{k-2}).
  $$
  Let $P=\sum_{k=0}^\infty p_kt^k$.  Then
  \begin{eqnarray*}
    P &=& 3t^m + 3t^{m+1} + t^mQ + 3\left(t^2+\cdots+t^{m-1}\right)Q\\
    &=& 3t^m\frac{1 - 2t + t^{m+1}}{1-3t+2t^m},
  \end{eqnarray*}
  again provided $|t|<\si_m$.

  The formula for the Hausdorff dimension of the infinite IFS is just
  $s=\log\tau_m/\log\om_m$, where by \cite[Corollary~3.17]{MU},
  $\tau_m$ is the supremum of all $x$ such that $\sum_k p_kx^k<1$,
  i.e., the (unique) positive root of $\sum_k p_kx^k=1$.  In other
  words, it is the unique solution in $(0,\si_m)$ of
  $$
  3\tau^m\frac{1 - 2\tau + \tau^{m+1}}{1 - 3\tau + \tau^m}=1.
  $$
  Rearranging this equation, one finds
  $$
  (3\tau^{m+1}-3\tau+1)\mathfrak{C}=0,
  $$
  where $\mathfrak{C}$ is the polynomial $\mathfrak{C} =
  1+t+\cdots+t^m$, none of whose roots are positive. It follows that
  $$
  \dim_H(\A_{\om_m}) = \log\tau_m/\log\om_m.
  $$
  It remains to show that $\dim_H(\A_{\om_m}) =
  \dim_H(\S_{\om_m})$.  The argument is similar to that for $\om_2$:
  it suffices to evaluate the growth of the number of hexagons, which
  follows from the generating function $Q$ found above. Indeed, the
  growth of the coefficient is asymptotically $h_k\asymp \si_m^{-k}$
  since $\si_m$ is the smallest root of the denominator of $Q$ (the
  radius of convergence mentioned above).  Thus, by Lemma
  \ref{lem-auxapp},
  $$
  \dim_H(\U_{\om_m}) = \frac{\log \si_m}{\log\om_m} < \frac{\log
    \tau_m}{\log\om_m} = \dim(\A_{\om_m}).
  $$
  The argument showing that $\U_{\om_m}$ is indeed the set of
  uniqueness, is analogous to the case $m=2$, so we omit it.
  Theorems~\ref{thm-hd} and \ref{thm:dimU} are now established.
\end{proof}

Thus, we have shown that ``almost every" point of $\S_{\om_m}$ (in
the sense of prevailing dimension) has at least {\em two}
different addresses. It is easy to prove a stronger claim:

\begin{prop}\label{prop:contin}
Define $\C_\la$ as the set of points in $\S_\la$ with less than a
continuum addresses, i.e.,
\[
\C_\la:=\left\{\mathbf x\in\S_\la :
\mathrm{card}\{\bm\e\in\Si^\infty : \mathbf
x=\xe\}<2^{\aleph_0}\right\}.
\]
Then
\[
\dim_H(\C_{\om_m})=\dim_H(\U_{\om_m}) < \dim_H(\S_{\om_m}).
\]
\end{prop}
\begin{proof}Let $\mathbf x=\xe$; if there exist an infinite
number of $k$'s such that $\e_k=i_k,\e_{k+1}=\dots=\e_{k+m}=j_k$
with $i_k\neq j_k$, then, obviously, $\mathbf x$ has a continuum
of addresses, because one can replace each $i_kj_k^m$ by
$j_ki_k^m$ independently of the rest of the address.

Thus, $\mathbf x=\xe\in\C_{\om_m}$ only if the tail of $\bm\e$ is
either $j^\infty$ for some $j\in\Si$ (a countable set we discard)
or a sequence $\bm\e'$ such that $\mathbf
x_{\bm\e'}\in\U'_{\om_m}$. Hence $\C_{\om_m}$ contains a
(countable) union of images of $\U'_{\om_m}$, each having the same
Hausdorff dimension, whence
$\dim_H(\C_{\om_m})=\dim_H(\U_{\om_m})$.
\end{proof}

We conjecture that the same claim is true for each
$\la\in(1/2,1)$. For the one-dimensional model this was shown by
the third author \cite{Sid}. Note that combinatorial questions of
such a kind make sense for $\la\ge2/3$ as well, since here the
holes are unimportant.


\bigskip

\begin{minipage}[t]{0.4\textwidth}
\obeylines
\small
\textsl{
  Department of Mathematics,
  UMIST,
  P.O. Box 88,
  Manchester M60 1QD,
  United Kingdom.
}
\end{minipage}
\hfill
\begin{minipage}[t]{0.4\textwidth}
\obeylines
\small
\textsl{
  E-mails:
  d.s.broomhead@umist.ac.uk
  j.montaldi@umist.ac.uk
  nikita.a.sidorov@umist.ac.uk
}
\end{minipage}

\end{document}